\documentclass{amsart}

\usepackage{amsmath,amsfonts,amssymb,amsthm,graphicx,cite,mathrsfs}

\newcommand{\lp}{\left(}
\newcommand{\rp}{\right)}

\newtheorem{theorem}{Theorem}[section]
\newtheorem{lemma}{Lemma}[section]
\newtheorem{proposition}{Proposition}[section]
\newtheorem{corollary}{Corollary}[section]
\newtheorem{definition}{Definition}[section]

\theoremstyle{remark}
\newtheorem{remark}{Remark}[section]

\begin{document}

\title[A basis for the eigenfunctions of deformed CMS operators]{A
basis for the polynomial eigenfunctions of deformed
Calogero-Moser-Sutherland operators}

\author{Martin Halln\"as} \address{SISSA, Via Beirut 2-4, 34014
Trieste, Italy} \email{hallnas@sissa.it}

\date{\today}

\begin{abstract}
We construct a linear basis for the polynomial eigenfunctions of a
family of deformed Calogero-Moser-Sutherland operators naturally
associated with hypergeometric polynomials. In our
construction the eigenfunctions are obtained as linear combinations
of polynomials which generalise the (super) Schur polynomials.  As a
byproduct, we obtain explicit series representations for the super
Jack polynomials.
\end{abstract}

\maketitle

\section{Introduction}
The main purpose of this paper is to construct and study a particular
linear basis for the polynomial eigenfunctions of a certain family of
partial differential operators which, following Sergeev and Veselov
\cite{SergVes2}, we will refer to as deformed
Calogero-Moser-Sutherland (CMS) operators. In our construction the
eigenfunctions are expressed as linear combinations of particular
polynomials which generalise the so-called (super) Schur polynomials;
see e.g.\ Fulton and Pragacz \cite{FultPrag}. Our main motivation for
the construction is that it leads to rather simple and explicit
formulae; see Section 7 for concrete examples.

To give a precise definition of the family of deformed CMS operators
we consider, and to better describe our results, we start by
discussing a simple and well known result from the theory of
polynomials in one variable: suppose that we are given a sequence of
polynomials
\begin{equation*}
  p_0(x), p_1(x),\ldots,p_n(x),\ldots,
\end{equation*}
where each polynomial $p_n$ is such that it has precisely degree $n$
and is an eigenfunction of a second order ordinary differential
operator
\begin{equation*}
  \mathcal{L} = \alpha(x)\frac{\partial^2}{\partial x^2} +
  \beta(x)\frac{\partial}{\partial x}
\end{equation*}
for some fixed polynomials $\alpha$ and $\beta$. It is then a
straightforward exercise to verify that $\alpha$ is of at most degree
two and that $\beta$ is of at most degree one, i.e.,
\begin{equation*}
  \alpha(x) = \alpha_2x^2 + \alpha_1x + \alpha_0,\quad \beta(x)
  = \beta_1x + \beta_0
\end{equation*}
for some (real) coefficients $\alpha_k$ and $\beta_\ell$; see e.g.\
Bochner \cite{Bochner}. Examples of such sequences of polynomials are
given by the classical orthogonal Hermite-, Laguerre- and Jacobi
polynomials, as well as the Bessel polynomials, all of which can be expressed in terms of hypergeometric functions. For a comprehensive
discussion of the classical orthogonal polynomials see for example
Andrews et al.\ \cite{AAR}, and for the Bessel polynomials the book by
Grosswald \cite{Gross}. It is interesting to note that, as long as
$\alpha$ is not identically zero, we can always reduce to one of these
four cases by an affine transformation of the variable $x$; see e.g.\
Bochner (\emph{loc.\ cit.}).

This type of (complete) sequences of polynomials have a natural
many-variable generalisation within the theory of symmetric
polynomials. In fact, Lassalle \cite{Lass1,Lass2,Lass3} and Macdonald
\cite{Macd1} introduced and studied a many-variable generalisations of
the classical orthogonal Hermite-, Laguerre- and Jacobi- polynomials
as eigenfunctions of partial differential operators
\begin{equation*}
  \mathcal{L}_n = \sum_{k=0}^2\alpha_kD^k_n +
  \sum_{\ell=0}^1\beta_\ell E^\ell_n
\end{equation*}
which can be obtained from the corresponding ordinary differential
operators $\mathcal{L}$ by replacing each term $x^k\partial^2/\partial
x^2$ by
\begin{equation*}
  D^k_n = \sum_{i=1}^nx_i^k\frac{\partial^2}{\partial x_i^2} +
  2\theta\sum_{i\neq j}\frac{x_i^k}{x_i - x_j}\frac{\partial}{\partial
  x_i}
\end{equation*}
and each term $x^\ell \partial/\partial x$ by
\begin{equation*}
  E^\ell_n = \sum_{i=1}^nx^\ell\frac{\partial}{\partial x_i}
\end{equation*}
for $k = 0,1,2$ and $\ell = 0,1$, respectively. These many-variable
polynomials have subsequently been extensively studied in the
literature. In particular, by Baker and Forrester \cite{BakerForr} and
van Diejen \cite{vanDie}. We also mention Heckman and Opdam's closely
related root system generalisation of the Jacobi polynomials; see
e.g.\ their paper \cite{HeckOp}. For additional related references see
e.g.\ the book by Dunkl and Xu \cite{DunklXu}.

In this paper we consider a further natural generalisation of the
ordinary differential operators $\mathcal{L}$ in two sequences of
independent variables $x = (x_1,\ldots,x_n)$ and $\tilde{x} =
(\tilde{x}_1,\ldots,\tilde{x}_{\tilde{n}})$. More precisely, we
consider the partial differential operators
\begin{equation}\label{CMSOps}
  \mathcal{L}_{n,\tilde{n}} = \sum_{k=0}^2 \alpha_k D^k_{n,\tilde{n}}
  + \sum_{\ell = 0}^1 \beta_\ell E^\ell_{n,\tilde{n}}
\end{equation}
obtained from $\mathcal{L}$ by replacing each term $x^k
\partial^2/\partial x^2$ by
\begin{align*}
  D^k_{n,\tilde{n}} &= \sum_{i=1}^nx^k_i\frac{\partial^2}{\partial
  x_j^2} - \theta\sum_{I=1}^{\tilde n}\tilde
  x^k_I\frac{\partial^2}{\partial\tilde x_I}\\ &\quad +
  2\theta\sum_{i\neq j}\frac{x_i^k}{x_i - x_j}\frac{\partial}{\partial
  x_i} - 2\sum_{I\neq J}\frac{\tilde{x}_I^k}{\tilde{x}_I -
  \tilde{x}_J}\frac{\partial}{\partial\tilde{x}_I}\\ &\quad -
  2\sum_{i,I}\frac{1}{x_i - \tilde x_I}\lp
  x_i^k\frac{\partial}{\partial x_i} + \theta\tilde
  x_I^k\frac{\partial}{\partial\tilde x_I}\rp + k(1 +
  \theta)\sum_{i=1}^nx_i^{k-1}\frac{\partial}{\partial x_i},
\end{align*}
and each term $x^\ell\partial/\partial x$ by
\begin{equation*}
  E^\ell_{n,\tilde{n}} = \sum_{i=1}^nx_i^\ell\frac{\partial}{\partial
  x_i} + \sum_{I=1}^{\tilde n}\tilde
  x_I^\ell\frac{\partial}{\partial\tilde x_I}
\end{equation*}
for $k = 0,1,2$ and $\ell = 0,1$, respectively. This type of partial
differential operators were introduced and studied by Chalykh et al.\
\cite{ChaFeiVes1} for $\tilde{n} = 1$ and by Sergeev
\cite{Sergeev1,Sergeev2} for arbitrary $\tilde{n}$. In special cases
they have subsequently been extensively studied by Sergeev and Veselov
\cite{SergVes1,SergVes2}. As mentioned above, following Sergeev and
Veselov we will refer to these partial differential operators
$\mathcal{L}_{n,\tilde{n}}$ as deformed Calogero-Moser-Sutherland
(CMS) operators. We mention that the operators studied in the papers
by Sergeev, as well as in the former of the two papers by Sergeev and
Veselov, are constructed using certain deformations of so-called
generalised root systems; see Serganova \cite{Serganova}. For the
classical series of such deformed generalised root systems the
resulting deformed CMS operators can also be obtained by specialising
to particular polynomials $\alpha$ and $\beta$ in
$\mathcal{L}_{n,\tilde{n}}$; see Table \ref{relationTable} below and
the paper \cite{HallLang} for specific examples. However, it is
interesting to note that the operator $\mathcal{L}_{n,\tilde{n}}$
corresponding to the Bessel polynomials (see Table
\ref{relationTable}) can not be directly defined in terms of root
systems.

\begin{table}[htbp]\label{relationTable}
  \begin{center}
    \begin{tabular}{|c|c|c|c|c|}
\hline
& & & &\\[-1.0ex]

$\alpha(x)$ & $\beta(x)$ & Eigenfuncs.\ of~$\mathcal{L}$ &
Type of & Root system\\
& & & potential &\\[1.2ex] \hline

& & & &\\[-1.0ex]

$x^2$ & $0$ & $x^n$ & Trigonometric & $A(n-1,\tilde{n}-1)$\\
& & (Monomials) & &\\[1.2ex]

$1$ & $-2x$ & $H_n(x)$ & Rational & $A(n-1,\tilde{n}-1)$\\
& & (Hermite pols.) & &\\[1.2ex]

$x$ & $a + 1 - x$ & $L^{(a)}_n(x)$ & Rational & $B(n,\tilde{n})$\\
& & (Laguerre pols.) & &\\[1.2ex]

$(1 - x^2)$ & $b - a -$ & $P^{(a,b)}_n(x)$ &
Trigonometric & $BC(n,\tilde{n})$
\\ & $(a + b + 2)x$ & (Jacobi pols.) & &\\[1.2ex]

$x^2$ & $b + ax$ & $y^{(a,b)}_n(x)$ & -- & -- \\
& & (gen. Bessel pols.) & &\\[1.2ex]

\hline
    \end{tabular}

\bigskip

\caption{Special cases of the deformed CMS operators
$\mathcal{L}_{n,\tilde{n}}$. In each case, the two rightmost columns
refers to the type of potential and root system in the construction
used by Sergeev and Veselov (\emph{loc.\ cit.}).}
  \end{center}
\end{table}

As indicated above, the main purpose of this paper is to construct and
study a particular linear basis for the polynomial eigenfunctions of
each of the deformed CMS operators $\mathcal{L}_{n,\tilde{n}}$. The
eigenfunctions are in this construction expressed as linear
combinations of certain polynomials generalising, or more accurately,
deforming the (super) Schur polynomials; see e.g.\ Fulton and Pragacz
\cite{FultPrag} for a discussion of the super Schur polynomials. Our
main motivation for using these polynomials is that they lead to
rather simple and explicit formulae. The type of series
representations we construct were first obtained by Langmann
\cite{Lang} (see also \cite{Lang2}) for eigenfunctions of an operator
$\mathcal{L}_n$ with only $\alpha_2$ and $\beta_1$ non-zero. His
results have subsequently been generalised to all `ordinary' CMS
operators $\mathcal{L}_n$, as well as their deformed counterparts
$\mathcal{L}_{n,\tilde{n}}$; see the paper \cite{HallLang} and
references therein. The present work is in many ways a natural
continuation of this latter paper. In particular, we obtain complete
proofs of a number of results which are only sketched or mentioned in
\cite{HallLang}. On the other hand, certain of the results obtained in
the present paper can be inferred from results in
\cite{HallLang}. However, the point of view in this latter paper, and
also in the papers by Langmann (\emph{loc.\ cit.}), is that of quantum
many-body systems of (deformed) Calogero-Sutherland type; see e.g.\
Calogero \cite{Calogero} and Sutherland \cite{Suth}. In all cases were
such an overlap of results occur our approach thus provides an
alternative and independent derivation, making no reference to their
relation with such quantum many-body systems.

We conclude this introduction by presenting a brief outline of the
paper. The first two sections are of an introductory nature. We begin
in Section 2 by giving a brief review of some basic facts from the
theory of symmetric functions. We also recall a definition of the Jack
polynomials and prove a short technical lemma related to the complete
symmetric polynomials. In Section 3 we recall the definition of the
so-called super Jack polynomials and review some of their basic
properties. In addition, we establish a particular triangular
structure in the expansion of the super Jack polynomials in ordinary
monomials. The main results of the paper are obtained in Sections
4-6. In Section 4 we establish particular identities which relates
each of the deformed CMS operators above to their adjoints for a
scalar product naturally associated with the super Jack
polynomials. In Section 5 we define and study certain polynomials
$f^{(m,\tilde{m})}_a(x,\tilde{x};\theta)$ labeled by two non-negative
integers $(m,\tilde{m})$, an integer vector
$a\in\mathbb{Z}^{m+\tilde{m}}$, and the parameter $\theta$. We prove
that they essentially coincide with the (super) Schur polynomials for
$\theta = 1$ and $\tilde{m} = 0$. We also study the structure of their
expansion in terms of super Jack polynomials. As a consequence, we
obtain a simple characterisation of the linear span of a certain
natural subset of the polynomials
$f^{(m,\tilde{m})}_a(x,\tilde{x};\theta)$. In Section 6 we construct
and study polynomial eigenfunctions of the deformed CMS operators
$\mathcal{L}_{n,\tilde{n}}$ as linear combinations of the polynomials
$f^{(m,\tilde{m})}_a(x,\tilde{x};\theta)$. An important aspect of our
construction is that the resulting eigenfunctions can be normalised
such that they are independent of the choice of parameters
$(m,\tilde{m})$. As we then discuss, this freedom in choosing the
values of $(m,\tilde{m})$ can be used to minimise the complexity of
the series representation of a given eigenfunction, in many cases
significantly below that of the canonical choice $(m,\tilde{m}) =
(n,\tilde{n})$. We also obtain a simple characterisation of the linear
span of the eigenfunctions we construct. In Section 7 we deduce the
explicit series expansion of the super Jack polynomials in terms of
the polynomials $f^{(\bar{m})}_a(x,\tilde{x};\theta)$. We also study
in some detail certain particularly simple special cases of this
series expansion. We conclude the paper in Section 8 by a brief
discussion of some open problems.

\section{Symmetric functions and Jack polynomials}
In this section we briefly recall some basic facts and definitions
concerning symmetric functions and Jack polynomials. We also prove a
short technical lemma on the so-called complete symmetric
polynomials. With a few minor exceptions we follow the notation of
Macdonald \cite{Macd2} to which the reader is referred for further
details.

Consider the algebra $\mathbb{C}\lbrack x_1,\ldots,x_n\rbrack$ of
polynomials in $n$ independent variables $x = (x_1,\ldots,x_n)$ with
complex coefficients. The subalgebra of all symmetric polynomials is
denoted $\Lambda_n$. It is graded by the degree of the polynomials,
i.e.,
\begin{equation*}
  \Lambda_n = \bigoplus_{k\geq 0}\Lambda^k_n
\end{equation*}
with $\Lambda^k_n$ the homogeneous component of $\Lambda_n$ of degree
$k$. Let $n\geq m$ and consider the homomorphism
\begin{equation*}
  \mathbb{C}\lbrack x_1,\ldots,x_n\rbrack\rightarrow\mathbb{C}\lbrack
  x_1,\ldots,x_m\rbrack
\end{equation*}
which sends each of the variables $x_{m+1},\ldots,x_{n}$ to zero and
the remaining variables $x_i$ to themselves. Let $\rho_{n,m}$ and
$\rho^k_{n,m}$ denote the restriction of this homomorphism to
$\Lambda_n$ and $\Lambda^k_n$, respectively. The inverse limit
\begin{equation*}
  \Lambda^k = \varprojlim\Lambda^k_n
\end{equation*}
of the linear spaces $\Lambda^k_n$ relative to the homomorphisms
$\rho^k_{m,n}$ can now be formed, and the algebra of symmetric
functions can be defined as the direct sum
\begin{equation*}
  \Lambda = \bigoplus_{k\geq 0}\Lambda^k.
\end{equation*}

A partition $\lambda = (\lambda_1,\lambda_2,\ldots)$ is any sequence
of non-negative integers in decreasing order, i.e.,
\begin{equation*}
  \lambda_1\geq\lambda_2\geq\cdots\geq\lambda_i\geq\cdots,
\end{equation*}
containing only a finite number of non-zero terms. These non-zero
terms $\lambda_i$ are called the parts of $\lambda$ and the number of
parts the length of $\lambda$, in the following denoted
$\ell(\lambda)$. The sum $|\lambda| := \lambda_1 + \lambda_2 +\cdots$
of its parts is referred to as the weight of $\lambda$. If $|\lambda|
= n$ it is said that $\lambda$ is a partition of $n$. We will for
simplicity not distinguish two partitions differing only by a string
of zeros at the end. A partition $\lambda$ can be identified with its
diagram, which consists of the points $(i,j)\in\mathbb{Z}^2$ such that
$1\leq j\leq\lambda_i$. The partition $\lambda^\prime$, obtained by
reflection in the main diagonal, is called the conjugate of
$\lambda$. On the set of the partitions of a given non-negative
integer $n$ is the so-called dominance order defined by
\begin{equation*}
  \mu\leq\lambda\Leftrightarrow \mu_1 + \cdots + \mu_i\leq \lambda_1 +
  \cdots + \lambda_i,\quad \forall i\geq 1.
\end{equation*}
More generally, we will write $a\leq b$ for any two integer vectors $a
= (a_1,\ldots,a_n)$ and $b = (b_1,\ldots,b_n)$ such that
\begin{equation}\label{intVectDominanceOrder}
  a_1 +\cdots + a_i\leq b_1 +\cdots + b_i
\end{equation}
for all $i = 1,\ldots, n$

Throughout the paper we will write $x^a = x_1^{a_1}\cdots x_n^{a_n}$
for any integer vector $a = (a_1,\ldots,a_n)\in\mathbb{Z}^n$. With
this notation in mind we recall the definition of the following linear
basises for $\Lambda_n$:
\begin{enumerate}
\item Monomial symmetric polynomials: defined for each partition
$\lambda$ of length $\ell(\lambda)\leq n$ by
\begin{equation*}
  m_\lambda(x_1,\ldots,x_n) = \sum_\alpha x^\alpha
\end{equation*}
where the sum extends over all distinct permutations $\alpha$ of
$\lambda$.
\item The power sums: defined for each $r\geq 1$ by
\begin{equation*}
  p_r(x_1,\ldots,x_n) = x_1^r + \ldots + x_n^r,
\end{equation*}
and for all partitions $\lambda$ such that $\ell(\lambda^\prime)\leq
n$,
\begin{equation*}
  p_\lambda = p_{\lambda_1}p_{\lambda_2}\cdots.
\end{equation*}
\item Elementary symmetric polynomials: can be defined for all $r\leq
  n$ by the expansion
\begin{equation*}
	\prod_{i=1}^n(1 + x_it) = \sum_{r=1}^n e_r(x_1,\ldots,x_n)t^r.
\end{equation*}
For all partitions $\lambda$ such that $\ell(\lambda^\prime)\leq n$,
\begin{equation*}
	e_\lambda = e_{\lambda_1}e_{\lambda_2}\cdots.
\end{equation*}
\item Complete symmetric polynomials: can be defined for all non-negative
integers $r$ by the expansion
\begin{equation*}
  \prod_{i=1}^n(1 - x_it)^{-1} = \sum_{r\geq 0}h_r(x_1,\ldots,x_n)t^r.
\end{equation*}
Furthermore, for each partition $\lambda$ such that
$\ell(\lambda)\leq n$,
\begin{equation*}
  h_\lambda = h_{\lambda_1}h_{\lambda_2}\cdots.
\end{equation*}
\item 'Modified' complete symmetric polynomials: can be defined for
  all non-negative integers $r$ and each real number $\theta$ by the
  expansion
\begin{equation*}
  \prod_{i=1}^n(1 - x_it)^{-\theta} = \sum_{r\geq 0}g_r(x_1,\ldots,x_n;\theta)t^r.
\end{equation*}
In addition, for all partitions $\lambda$ such that
$\ell(\lambda)\leq n$,
\begin{equation*}
  g_\lambda = g_{\lambda_1}g_{\lambda_2}\cdots.
\end{equation*}
\end{enumerate}

It is a well known fact that these polynomials
$m_\lambda(x_1,\ldots,x_n)$, $p_\lambda(x_1,\ldots,x_n)$,
$e_\lambda(x_1,\ldots,x_n)$, $h_\lambda(x_1,\ldots,x_n)$ and
$g_\lambda(x_1,\ldots,x_n)$, under the restrictions on the partitions
$\lambda$ stated above, all form linear basises for $\Lambda_n$; see
e.g.\ Sections I.2 and VI.10 in Macdonald \cite{Macd2}. In addition,
since these polynomials are stable under the homomorphisms
$\rho_{n,m}$, the symmetric functions $m_\lambda$, $p_\lambda$,
$e_\lambda$, $h_\lambda$ and $g_\lambda$ can be defined and form
linear basises for $\Lambda$.

In later parts of the paper we make use of the fact that the complete
symmetric polynomials in three variables can be expressed as a very
particular quotient. This fact is established in the following:

\begin{lemma}\label{threeLemma}
Set $h_r = 0$ for all $r<0$. Let $x$, $y$ and $z$ be three independent
variables. Then, for each non-negative integer $k$,
\begin{equation*}
  \frac{x^k(z - y) + y^k(x - z) + z^k(y - x)}{(y - x)(x - z)(z - y)} =
  -h_{k-2}(x,y,z).
\end{equation*}
\end{lemma}

\begin{proof}
Let
\begin{equation*}
  q_k(x,y,z) = \frac{x^k(z - y) + y^k(x - z) + z^k(y - x)}{(y - x)(x -
  z)(z - y)}
\end{equation*}
and observe that
\begin{equation*}
\begin{split}
  \sum_{k\geq 0}q_k(x,y,z)t^k &= \frac{1}{(y - x)(x - z)}\sum_{k\geq
  0}(xt)^k + \frac{1}{(y - x)(z - y)}\sum_{k\geq 0}(yt)^k\\ &\quad +
  \frac{1}{(x - z)(z - y)}\sum_{k\geq 0}(zt)^k\\ &= \frac{1}{(y - x)(x
  - z)(1 - xt)} + \frac{1}{(y - x)(z - y)(1 - yt)}\\ &\quad +
  \frac{1}{(x - z)(z - y)(1 - zt)}.
\end{split}
\end{equation*}
This latter sum can be written as a single fraction with nominator
\begin{multline*}
  (z - y)(1 - yt)(1 - zt) + (x - z)(1 - xt)(1 - zt) + (y - x)(1 -
  xt)(1 - yt)\\ = - t^2(y - x)(x - z)(z - y).
\end{multline*}
It follows that
\begin{equation*}
\begin{split}
  \sum_{k\geq 0}q_k(x,y,z)t^k &= -t^2(1 - xt)^{-1}(1 - yt)^{-1}(1 -
  zt)^{-1} \\ &= -\sum_{k\geq 0}h_k(x,y,z)t^{k+2}.
\end{split}
\end{equation*}
\end{proof}

We proceed to recall a definition of the Jack polynomials suitable for
our purposes. In doing so we use the inverse $\theta = 1/\alpha$ of
the parameter $\alpha$ used by Macdonald \cite{Macd2}. We recall that
for each `square' $s = (i,j)$ in the diagram of a partition $\lambda$
the so-called arm-length $a(s)$ and leg-length $l(s)$ are given by
$a(s) = \lambda_i - j$ and $l(s) = \lambda_j^\prime - i$,
respectively. We let
\begin{equation*}
  b_\lambda(\theta) = \prod_{s\in\lambda}\frac{a(s) + \theta l(s) +
  \theta}{a(s) + \theta l(s) + 1}.
\end{equation*}
Following Macdonald (\emph{loc.\ cit.}) we let $\Box^{1/\theta}_n$
be the so-called Laplace-Beltrami operator given by
\begin{equation*}
  \Box^{1/\theta}_n = \frac{1}{2\theta}D^2_n - (n - 1)E^1_n =
  \frac{1}{2\theta}\sum_ix_i^2\frac{\partial^2}{\partial x_i^2} +
  \sum_{i\neq j}\frac{x_ix_j}{x_i - x_j}\frac{\partial}{\partial x_i}.
\end{equation*}
It is clear that $\Box^{1/\theta}_n$ is stable under the restriction
homomorphisms $\rho_{n,m}$, i.e., that
\begin{equation*}
  \rho_{n,m}\circ\Box^{1/\theta}_n = \Box^{1/\theta}_m\circ\rho_{n,m}
\end{equation*}
for all $m\leq n$. It follows that its inverse limit
\begin{equation*}
  \Box^{1/\theta} = \varprojlim\Box^{1/\theta}_n
\end{equation*}
is a well defined operator on $\Lambda$. Jack's symmetric functions
$P_\lambda$ can now be defined by the following result due to
Macdonald (see e.g.\ Example 3 in Section VI.5 of \cite{Macd2}):

\begin{theorem}[Macdonald]\label{jackDef}
If $\theta$ is not a negative rational number or zero there exist for
each partition $\lambda$ a unique eigenfunction $P_\lambda$ of the
operator $\Box^{1/\theta}$ such that
\begin{equation*}
  P_\lambda = m_\lambda + \sum_{\mu<\lambda}u_{\lambda\mu}m_\mu
\end{equation*}
for some coefficients $u_{\lambda\mu}$.
\end{theorem}

By setting all variables $x_i = 0$ for $i>n$, where $n$ is some
positive integer, we obtain the corresponding Jack polynomials
$P_\lambda(x_1,\ldots,x_n)$ in $n$ variables $x = (x_1,\ldots,x_n)$.

Let $x = (x_1,x_2,\ldots)$ and $y = (y_1,y_2,\ldots)$ be two infinite
sequences of independent variables. A fundamental object in the theory
of Jack's symmetric functions is the infinite product
\begin{equation}\label{infiniteProd}
  \Pi(x,y;\theta) = \prod_{i,j}(1 - x_iy_j)^{-\theta}.
\end{equation}
We let $\Pi_{n,m}(x,y;\theta)$ denote its restriction to a finite
number of variables $x = (x_1,\ldots,x_n)$ and $y =
(y_,\ldots,y_m)$. The following well-known result is due to Stanley
\cite{Stan}:

\begin{proposition}[Stanley]\label{jackExpProposition}
  If $\theta$ is not a negative rational number or zero, then the
  infinite product $\Pi(x,y;\theta)$ has the following expansion in
  terms of Jack's symmetric functions:
\begin{equation}\label{JackExpansion}
  \Pi(x,y;\theta) = \sum_\lambda
  b_\lambda(\theta)P_\lambda(x;\theta)P_\lambda(y;\theta)
\end{equation}
where the sum extends over all partitions.
\end{proposition}

\section{The super Jack polynomials}
In this section we recall the definition of the so-called super Jack
polynomials, as stated by Kerov et al.\ \cite{KOO}. We also recall
certain related results due to Sergeev and Veselov
\cite{SergVes1,SergVes2}. In addition, we establish a particular
triangular structure in the expansion of the super Jack polynomials in
ordinary monomials, and we derive the analogue of Proposition
\ref{jackExpProposition} for the super Jack polynomials.

Let $x = (x_1,\ldots,x_n)$ and $\tilde{x} =
(\tilde{x}_1,\ldots,\tilde{x}_{\tilde{n}})$ be two sequences of
independent variables. Following Sergeev and Veselov \cite{SergVes1}
we let $\Lambda_{n,\tilde{n},\theta}$ be the subalgebra of
$\mathbb{C}\lbrack
x_1,\ldots,x_n,\tilde{x}_1,\ldots,\tilde{x}_{\tilde{n}}\rbrack$
consisting of all polynomials $p(x,\tilde{x})$ which, in addition to
being separately symmetric in the variables $x$ and $\tilde{x}$,
satisfy the condition
\begin{equation*}
  \left(\frac{\partial}{\partial x_i} + \theta\frac{\partial}{\partial
  \tilde{x}_I}\right)p(x,\tilde{x}) = 0
\end{equation*}
on each hyperplane $x_i = \tilde{x}_I$ with $i = 1,\ldots,n$ and $I =
1,\ldots,\tilde{n}$. In addition, for each non-negative integer $k$ we
let $\Lambda_{n,\tilde{n},\theta}^k$ be the homogeneous component of
$\Lambda_{n,\tilde{n},\theta}$ of degree $k$. It is clear that the
`deformed' power sums
\begin{equation*}
  p_{r,\theta}(x,\tilde{x}) = x_1^r + \cdots + x_n^r -
  \theta^{-1}\left(\tilde{x}_1^r + \cdots +
  \tilde{x}_{\tilde{n}}^2\right)
\end{equation*}
for $r\geq 1$ are contained in $\Lambda_{n,\tilde{n},\theta}$. We let
$\varphi_{n,\tilde{n}}:
\Lambda\rightarrow\Lambda_{n,\tilde{n},\theta}$ be the algebra
homomorphism defined by
\begin{equation*}
  \varphi_{n,\tilde{n}}(p_r) = p_{r,\theta}(x,\tilde{x})
\end{equation*}
for all $r\geq 1$. Note that since the the power sums $p_r$ are free
generators of $\Lambda$, this uniquely determines
$\varphi_{n,\tilde{n}}$. The super Jack polynomials can now be defined
as follows:

\begin{definition}[Kerov et al.]
For each partition $\lambda$ the super Jack polynomial
$SP_\lambda(x,\tilde{x})$ is defined by
\begin{equation*}
  SP_\lambda(x,\tilde{x}) = \varphi_{n,\tilde{n}}(P_\lambda).
\end{equation*}
\end{definition}

As a direct consequence of this definition and the fact that a Jack
polynomial $P_\lambda(x)$, if non-zero, is homogeneous of degree
$|\lambda|$ we obtain the following:

\begin{lemma}\label{homLemma}
If a super Jack polynomial $SP_\lambda$ is non-zero it is homogeneous
of degree $|\lambda|$.
\end{lemma}

We let $H_{n,\tilde{n}}$ be the set of partitions contained in the fat
$(n, \tilde{n})$-hook, i.e., the set of partitions $\lambda =
(\lambda_1,\lambda_2,\ldots)$ such that
$\lambda_{n+1}\leq\tilde{n}$. For any partition $\lambda\in
H_{n,\tilde{n}}$ we let
\begin{equation*}
  {}^n\lambda = (\lambda_1,\ldots,\lambda_n),\quad {}_n\lambda =
  (\lambda_{n+1},\lambda_{n+2},\ldots).
\end{equation*}
From Sergeev and Veselov \cite{SergVes2} we now recall the following:

\begin{theorem}[Sergeev and Veselov]
Assume that $\theta$ is not a negative rational number or zero. Then
the kernel of $\varphi_{n,\tilde{n}}$ is spanned by the Jack's
symmetric functions $P_\lambda$ indexed by the partitions
$\lambda\notin H_{n,\tilde{n}}$. Moreover, the super Jack
polynomials $SP_\lambda(x,\tilde{x})$, indexed by the partitions
$\lambda\in H_{n,\tilde{n}}$, form a linear basis for
$\Lambda_{n,\tilde{n},\theta}$.
\end{theorem}

We proceed to establish a particular triangular structure in the
expansion of the super Jack polynomials in ordinary monomials. This
result will play an important role in later parts of the paper.

\begin{lemma}\label{triangLemma}
Assume that $\theta$ is not a negative rational number or zero and
let $\lambda\in H_{\bar{n}}$. The super Jack polynomial
$SP_\lambda(x,\tilde{x};\theta)$ is then a linear combination of
monomials $x^a\tilde{x}^b$ with $(a,b)\leq
({}^n\lambda,{}_n\lambda^\prime)$. Moreover, the leading term is
given by
\begin{equation}\label{leadingSuperTerm}
  (-1)^{|{}_n\lambda|} b_{{}_n\lambda^\prime}\left(\theta^{-1}\right)
  x^{{}^n\lambda}\tilde{x}^{{}_n\lambda^\prime}.
\end{equation}
\end{lemma}

It has been shown by Sergeev and Veselov \cite{SergVes2} that
\eqref{leadingSuperTerm} is the leading term of $SP_\lambda$ in the
lexicographic order; see the proof of their Theorem 2. We recall that
for two integer vectors $a,b\in\mathbb{Z}^m$, where $m$ is some
positive integer, the integer vector $a$ is said to be of lower order
than $b$ in the lexicographic order if $a\neq b$ and the first
non-zero term in $b-a$ is positive. Note that this is the case if
$a<b$. Consequently, Lemma \ref{triangLemma} generalises the result of
Sergeev and Veselov. We will obtain a proof of the lemma by extending
their argument. The idea is to first deduce an explicit expansion of
the super Jack polynomials in terms of Jack- and so-called skew Jack
polynomials. The proof of Lemma \ref{triangLemma} is then obtained by
well-known properties of these latter polynomials.

\begin{proof}[Proof of Lemma \ref{triangLemma}]
In order for the proof to be self-contained we start by recalling
the definition of the skew Jack symmetric functions, as stated in
Section VI.7 of Macdonald \cite{Macd2}. To any Jack symmetric
function $P_\lambda(x;\theta)$ is associated
\begin{equation*}
  Q_\lambda(x;\theta) = b_\lambda(\theta)P_\lambda(x;\theta).
\end{equation*}
We define the scalar product $\langle\cdot,\cdot\rangle$ on $\Lambda$
(antilinear in its second argument) by setting
\begin{equation*}
  \langle P_\lambda,Q_\mu\rangle = \delta_{\lambda\mu}
\end{equation*}
for all partitions $\lambda$ and $\mu$. For any two partitions
$\lambda$ and $\mu$ the skew Jack symmetric functions
$P_{\lambda/\mu}$ can now be defined by requiring
\begin{equation*}
  \langle P_{\lambda/\mu},Q_\nu\rangle = \langle P_\lambda,Q_\mu
  Q_\nu\rangle
\end{equation*}
for all partitions $\nu$. The corresponding skew Jack polynomials
$P_{\lambda/\mu}(x_1,\ldots,x_n)$ are obtained by setting all
variables $x_i = 0$ for $i>n$.

We proceed to prove the statement. To this end, we let $x =
(x_1,x_2,\ldots)$ and $\tilde{x} = (\tilde{x}_1,\tilde{x}_2,\ldots)$
be two infinite sequences of independent variables. It is a well known
fact that
\begin{equation}\label{jackSkewExpansion}
  P_\lambda(x,\tilde{x}) =
  \sum_{\mu\subset\lambda}P_{\lambda/\mu}(x)P_\mu(\tilde{x});
\end{equation}
see e.g.\ Macdonald (\textit{loc.\ cit.}). We let $\omega_\theta$
denote the automorphism of $\Lambda$ defined by
\begin{equation*}
  \omega_\theta(p_r) = (-1)^{r-1}\theta p_r
\end{equation*}
for all $r\geq 1$, and recall from Macdonald (\emph{loc.\
cit.}) that
\begin{equation*}
  \omega_{\theta^{-1}}(P_\lambda(x;\theta)) =
  Q_{\lambda^\prime}(x;\theta^{-1}) = b_{\lambda^\prime}(\theta^{-1})
  P_{\lambda^\prime}(x;\theta^{-1}).
\end{equation*}
We observe that acting with the homomorphism $\varphi_{n,\tilde{n}}$
on a Jack's symmetric function $P_\lambda(x,\tilde{x})$ is equivalent
to first acting with the automorphism $\omega_{\theta^{-1}}$ in the
variables $\tilde{x}$, followed by changing the sign of all variables
$\tilde{x}$, and finally setting all variables $x_i = 0$ and
$\tilde{x}_I = 0$ for $i>n$ and $I>\tilde{n}$, respectively. It
follows that
\begin{equation}\label{superJackExpInJack}
  SP_\lambda(x,\tilde{x};\theta) =
  \sum_{\mu\subset\lambda}(-1)^{|\mu|}b_{\mu^\prime}(\theta^{-1})
  P_{\lambda/\mu}(x;\theta)P_{\mu^\prime}(\tilde{x};\theta^{-1})
\end{equation}
with $x = (x_1,\ldots,x_n)$ and $\tilde{x} =
(\tilde{x}_1,\ldots,\tilde{x}_{\tilde{n}})$. We recall that a skew
Jack polynomial $P_{\lambda/\mu}(x)$ has an expansion in `ordinary'
Jack polynomials of the form
\begin{equation*}
  P_{\lambda/\mu}(x) =
  \sum_{\nu\subset\lambda}c^\lambda_{\mu\nu}P_\nu(x),
\end{equation*}
where the coefficients $c^\lambda_{\mu\nu}$ are non-zero only if
$|\mu| + |\nu| = |\lambda|$; see e.g.\ Macdonald (\emph{loc.\
  cit.}). We insert this expansion in \eqref{superJackExpInJack} and
fix a non-zero term
$P_\nu(x;\theta)P_{\mu^\prime}(\tilde{x};\theta^{-1})$ appearing in
the resulting expression. Since $\nu\subset\lambda$ and $P_\nu(x)$
vanishes unless $\ell(\nu)\leq n$ we have that
$\nu\leq{}^n\lambda$. For any integer $a$ we let $\langle a\rangle =
\max(a,0)$. Since $P_{\lambda/\mu} = 0$ unless $0\leq\lambda^\prime_i
- \mu^\prime_i\leq n$ for all $i\geq 1$ (see e.g.\ Macdonald
(\emph{loc.\ cit.})), we have that
$\mu^\prime_i\geq\langle\lambda^\prime_i-n\rangle =
{}_n\lambda^\prime_i$ for all $i\geq 1$. It follows that
\begin{equation*}
\begin{split}
  \nu_1 + \cdots + \nu_n + \mu^\prime_1 + \cdots + \mu^\prime_i &=
  |\lambda| - \mu^\prime_{i+1} - \mu^\prime_{i+1}\cdots\\ &\leq
  |\lambda| - {}_n\lambda^\prime_{i+1} - {}_n\lambda^\prime_{i+2}
  -\cdots\\ &= {}^n\lambda_1 + \cdots + {}^n\lambda_n +
  {}_n\lambda^\prime_1 +\cdots + {}_n\lambda^\prime_i
\end{split}
\end{equation*}
for each $i\geq 1$. Theorem \ref{jackDef} and the fact that each
symmetric monomial $m_\mu(x)$ is a linear combination of monomials
$x^a$ with $a\leq\mu$ thus implies that the super Jack polynomial
$SP_\lambda(x,\tilde{x})$ indeed is a linear combination of monomials
$x^a\tilde{x}^b$ with $(a,b)\leq
({}^n\lambda,{}_n\lambda^\prime)$. There remains only to verify that
\eqref{leadingSuperTerm} is the leading term in the expansion of
$SP_\lambda(x,\tilde{x})$ in such monomials. However, using the fact
that $\mu_i^\prime\geq {}_n\lambda_i^\prime$ we infer from
\eqref{superJackExpInJack} and Theorem \ref{jackDef} that the monomial
$x^{{}^n\lambda}\tilde{x}^{{}_n\lambda^\prime}$ can appear only in the
term
$P_{\lambda/{}_n\lambda}(x)P_{{}_n\lambda^\prime}(\tilde{x})$. That it
indeed appears in this term, and with the coefficient given in the
statement, follows from \eqref{jackSkewExpansion}.
\end{proof}

With $y = (y_1,\ldots,y_m)$ and $\tilde{y} =
(\tilde{y}_1,\ldots,\tilde{y}_{\tilde{m}})$ two sequences of
independent variables we now obtain the analogue of Proposition
\ref{jackExpProposition} for the super Jack polynomials.

\begin{proposition}\label{superJackExpansionProposition}
Assume that $\theta$ is not a negative rational number or zero and let
\begin{equation*}
  \Pi_{n,\tilde{n},m,\tilde{m}}(x,\tilde{x},y,\tilde{y};\theta) =
  \frac{\prod_{i,J}(1 - x_i\tilde{y}_J)\prod_{I,j}(1 -
  \tilde{x}_Iy_j)}{\prod_{i,j}(1 - x_iy_j)^\theta\prod_{I,J}(1 -
  \tilde{x}_I\tilde{y}_J)^{1/\theta}}.
\end{equation*}
Then
\begin{equation*}
  \Pi_{n,\tilde{n},m,\tilde{m}}(x,\tilde{x},y,\tilde{y};\theta) =
  \sum_\lambda b_\lambda(\theta)
  SP_\lambda(x,\tilde{x})SP_\lambda(y,\tilde{y})
\end{equation*}
where the sum is over all partitions $\lambda\in H_{n,\tilde{n}}\cap
H_{m,\tilde{m}}$.
\end{proposition}

\begin{proof}
To prove the statement we compute the action of the homomorphism
$\varphi = \varphi_{n,\tilde{n}}$ on the infinite product
$\Pi(x,y;\theta)$ in both the variables $x$ and $y$.

In the proof of Lemma 3 in \cite{SergVes2} Sergeev and Veselov
established that
\begin{equation}\label{xAction}
  \varphi_x\left(\Pi(x,y;\theta)\right) = \prod_j\prod_{i=1}^n(1 -
  x_iy_j)^{-\theta}\prod_{I=1}^{\tilde{n}}(1 - \tilde{x}_Iy_j),
\end{equation}
where the suffix $x$ indicates that the homomorphism $\varphi$ acts in
the variables $x$. To deduce the effect of a subsequent application of
$\varphi_y$ we essentially repeat their computation. We first observe
that since $\varphi$ is a homomorphism, it is sufficient to compute
its action on
\begin{equation*}
  \prod_j(1 - x_iy_j)^{-\theta}\quad\text{and}\quad \prod_j(1 -
  \tilde{x}_Iy_j)
\end{equation*}
for some fixed $i$ and $I$. For the first product \eqref{xAction}
directly implies that
\begin{equation*}
  \varphi_y\left(\prod_j(1 - x_iy_j)^{-\theta}\right) =
  \prod_{j=1}^n(1 - x_iy_j)^{-\theta}\prod_{J=1}^{\tilde{n}}(1 -
  x_i\tilde{y}_J).
\end{equation*}
To treat the second product we let $\sigma_\theta$ be the automorphism
of $\Lambda$ defined by
\begin{equation*}
  \sigma_\theta(p_r) = -\theta^{-1}p_r
\end{equation*}
for all $r\geq 1$. With an infinite number of independent variables
$\tilde{y} = (\tilde{y}_1,\tilde{y}_2,\ldots)$ we thus obtain
\begin{equation*}
\begin{split}
  \varphi_y\left(\prod_j(1 - \tilde{x}_Iy_j)\right) &= \prod_{j=1}^n(1
  - \tilde{x}_Iy_j)\sigma_{\theta,\tilde{y}}\left(\prod_J(1 -
  \tilde{x}_I\tilde{y}_J)\right)\\ &= \prod_{j=1}^n(1 -
  \tilde{x}_Iy_j)\sigma_{\theta,\tilde{y}}\left(\exp\log\prod_J(1 -
  \tilde{x}_I\tilde{y}_J)\right)\\ &= \prod_{j=1}^n(1 -
  \tilde{x}_Iy_j)\sigma_{\theta,\tilde{y}}\left(\exp\sum_{r\geq
  1}\frac{\tilde{x}_I^rp_r(\tilde{y})}{r}\right)\\ &= \prod_{j=1}^n(1
  - \tilde{x}_Iy_j)\exp\left(-\theta^{-1}\sum_{r\geq
  1}\frac{\tilde{x}_I^rp_r(\tilde{y})}{r}\right)\\ &= \prod_{j=1}^n(1
  - \tilde{x}_Iy_j)\prod_J(1 - \tilde{x}_I\tilde{y}_J)^{-1/\theta}.
\end{split}
\end{equation*}
Setting all variables $\tilde{y}_J = 0$ for $J>\tilde{n}$ the proof
now follows from Proposition \ref{jackExpProposition} and the
definition of the super Jack polynomials.
\end{proof}

\section{Identities and adjoints for deformed CMS operators}
The main purpose of this section is to generalise the well known
identity
\begin{equation}\label{StanleysId}
  D^2_{n,x}\Pi_{n,n}(x,y) = D^2_{n,y}\Pi_{n,n}(x,y),
\end{equation}
where the suffix's $x$ and $y$ indicates that the operator acts in the
variables $x$ and $y$, respectively, to similar identities for each of
the operators $\mathcal{L}_{n,\tilde{n}}$. We mention that the
identity \eqref{StanleysId} can be obtained as a direct consequence of
the expansion of $\Pi_{n,n}$ in terms of Jack polynomials (stated in
Proposition \ref{jackExpProposition}) and the fact that the Jack
polynomials are eigenfunctions of $D^2_n$ (see e.g.\ Theorem 3.1 in
Stanley \cite{Stan}). In addition, it is a well known fact that the
identity \eqref{StanleysId} is equivalent to the self-adjointness of
$D^2_n$ for a particular scalar product on $\Lambda_n$ (for which the
Jack polynomials are pairwise orthogonal); see Statement 3.11 and
Example 3 in Section VI.3 of Macdonald \cite{Macd2}. In this section
we also establish a natural generalisation of this fact to the
operators $\mathcal{L}_{n,\tilde{n}}$.

To simplify the exposition in this and following sections we will at
this point change our notation somewhat. We set $x_{n+I} =
\tilde{x}_I$ for $I = 1,\ldots,\tilde{n}$, and let $p$ denote the
`parity' function defined by
\begin{equation*}
  p(i) = \begin{cases}
    0 & \text{if}~0\leq i\leq n\\
    1 & \text{if}~n<i\leq n+\tilde{n}.
  \end{cases}
\end{equation*}
We also collect the (non-negative) integers $n$ and $\tilde{n}$ in a
vector $\bar{n} = (n,\tilde{n})$, and similarly for $m$ and
$\tilde{m}$. In addition, we let $|\bar{n}| = n + \tilde{n}$ and
$|\bar{m}| = m + \tilde{m}$. The operators $E^\ell_{n,\tilde{n}}$ and
$D^k_{n,\tilde{n}}$ can now be rewritten in the following simple form:
\begin{equation*}
  E^\ell_{\bar{n}} =
  \sum_{i=1}^{|\bar{n}|}x_i^\ell\frac{\partial}{\partial x_i}
\end{equation*}
and
\begin{multline*}
  D^k_{\bar{n}} =
  \sum_{i=1}^{|\bar{n}|}(-\theta)^{p(i)}x_i^k\frac{\partial^2}{\partial
  x_i^2} - 2\sum_{i\neq
  j}(-\theta)^{1-p(j)}\frac{x_i^k}{x_i-x_j}\frac{\partial}{\partial
  x_i}\\ + k\sum_{i=1}^{|\bar{n}|}\lp 1 - (-\theta)^{1-p(i)}\rp
  x_i^{k-1}\frac{\partial}{\partial x_i}.
\end{multline*}
Similarly, setting $y_{m+J} = \tilde{y}_{J}$ for $J =
1,\ldots\tilde{m}$ and introducing a corresponding `parity' function
\begin{equation*}
  q(j) = \begin{cases}
    0 & \text{if}~0\leq j\leq m\\
    1 & \text{if}~m < j\leq m+\tilde{m},
  \end{cases}
\end{equation*}
the function
$\Pi_{n,\tilde{n},m,\tilde{m}}(x,\tilde{x},y,\tilde{y};\theta)$ can be
rewritten as
\begin{equation*}
  \Pi_{\bar n,\bar m}(x,y;\theta) = \prod_{i,j}(1 -
  x_iy_j)^{(-\theta)^{1-p(i)-q(j)}}.
\end{equation*}
We proceed to associate a particular operator in the variables $y$ to
each of the operators $E^\ell_{\bar{n}}$ and $D^k_{\bar{n}}$. In doing
so we make use of the short hand notation $n_\theta = \theta n -
\tilde{n}$ and $m_\theta = \theta m - \tilde{m}$. For $\ell = 0,1$ we
let
\begin{equation*}
  \bar{E}^\ell_{\bar{n},\bar{m}} =
  \sum_{j=1}^{|\bar{m}|}y_j^{2-\ell}\frac{\partial}{\partial y_j} + (1
  - \ell)n_\theta p_{1,\theta}(y),
\end{equation*}
and for $k = 0,1,2$ we let
\begin{multline*}
  \bar{D}^k_{\bar{n},\bar{m}} =
  \sum_{j=1}^{|\bar{m}|}(-\theta)^{q(j)}y_j^{4-k}
  \frac{\partial^2}{\partial y_j^2} - 2\sum_{j\neq
  l}(-\theta)^{1-q(l)}\frac{y_j^{4-k}}{y_j-y_l}\frac{\partial}{\partial
  y_j}\\ +
  \sum_{j=1}^{|\bar{m}|}C_{j,k}y_j^{3-k}\frac{\partial}{\partial y_j}
  + P_k(y)
\end{multline*}
with the constants
\begin{equation*}
  C_{j,k} = 2(n_\theta - m_\theta) + (2 - k)\lp(-\theta)^{q(j)} -
  (-\theta)^{1-q(j)}\rp + k\lp 1 - (-\theta)^{1-q(j)}\rp
\end{equation*}
and where the polynomials
\begin{equation*}
  P_k(y) = (1 - \delta_{k,2})n_\theta(n_\theta + 1) p_{2-k,\theta}(y)
  + \delta_{k,0}\theta n_\theta\lp p^2_{1,\theta}(y) -
  p_{2,\theta}(y)\rp.
\end{equation*}
We note that the polynomials $P_k$ are contained in
$\Lambda_{\bar{m},\theta}^{2-k}$. We are now ready to state and prove
identities which generalise \eqref{StanleysId} to the operators
$E^{\ell}_{\bar{n}}$ and $D^k_{\bar{n}}$.

\begin{lemma}\label{IdLemma}
For $\ell = 0,1$,
\begin{equation}\label{EId}
  E^\ell_{\bar{n},x}\Pi_{\bar{n},\bar{m}}(x,y) =
  \bar{E}^\ell_{\bar{n},\bar{m},y}\Pi_{\bar{n},\bar{m}}(x,y),
\end{equation}
and for $k = 0,1,2$,
\begin{equation}\label{DId}
  D^k_{\bar{n},x}\Pi_{\bar{n},\bar{m}}(x,y) =
  \bar{D}^k_{\bar{n},\bar{m},y}\Pi_{\bar{n},\bar{m}}(x,y).
\end{equation}
\end{lemma}

\begin{proof}
The specific number of variables $x$ and $y$ is in most parts of the
proof of little importance. For simplicity of notation, we therefore
suppress the index's $\bar{n}$ and $\bar{m}$. In addition, we assume
that, whenever they occur, $\ell = 0,1$ and $k = 0,1,2$.

We start by proving the identity \eqref{EId}. To this end we observe
that
\begin{equation}\label{firstDerivative}
  \Pi^{-1}(x,y)\frac{\partial}{\partial x_i}\Pi(x,y) =
  -\sum_j(-\theta)^{1-p(i)-q(j)}\frac{y_j}{1 - x_iy_j}
\end{equation}
and similarly for $\partial/\partial y_j$. This, together with the
fact that
\begin{equation*}
  \frac{x_iy_j^{2-\ell} - x_i^\ell y_j}{1 - x_iy_j} = (\ell - 1)y_j,
\end{equation*}
implies the equalities
\begin{equation*}
\begin{split}
  \Pi^{-1}(x,y)\lp E^\ell_x - E^{2-\ell}_y\rp\Pi(x,y) &=
  \sum_{i,j}(-\theta)^{1-p(i)-q(j)}\frac{x_iy_j^{2-\ell} - x_i^\ell
  y_j}{1 - x_iy_j}\\ &= (\ell -
  1)\sum_{i,j}(-\theta)^{1-p(i)-q(j)}y_j\\ &= (1 - \ell)(n\theta -
  \tilde n)p_{1,\theta}(y),
\end{split}
\end{equation*}
from which \eqref{EId} immediately follows.

We turn now to the proof of the identity \eqref{DId}. We begin with
the case $k = 2$ and observe that
\begin{multline*}
  \Pi^{-1}(x,y)(-\theta)^{p(i)}x_i^2\frac{\partial^2}{\partial
  x_i^2}\Pi(x,y)\\ =
  \sum_j(-\theta)^{1-q(j)}\lp(-\theta)^{1-p(i)-q(j)} -
  1\rp\frac{x_i^2y_j^2}{(1 - x_iy_j)^2}\\ + \sum_{j\neq
  l}(-\theta)^{2-p(i)-q(j)-q(l)}\frac{x_i^2y_jy_l}{(1 - x_iy_j)(1 -
  x_iy_l)}.
\end{multline*}
This, together with \eqref{firstDerivative} and the corresponding
results for $(-\theta)^{q(j)}y_j^2\partial^2/\partial y_j^2$ and
$\partial/\partial y_j$, imply
\begin{multline}\label{diffAction}
  \Pi^{-1}(x,y)\lp D^2_x - \bar{D}^2_y\rp\Pi(x,y)\\ =
  \sum_i\sum_{j\neq
  l}(-\theta)^{2-p(i)-q(j)-q(l)}\lp\frac{x_i^2y_jy_l}{(1 - x_iy_j)(1 -
  x_iy_l)} - 2\frac{x_iy_j^2}{(y_j - y_l)(1 - x_iy_j)}\rp\\ +
  \sum_{i\neq j}\sum_l(-\theta)^{2-p(i)-p(j)-q(l)}\lp
  2\frac{x_i^2y_l}{(x_i - x_j)(1 - x_iy_l)} - \frac{x_ix_jy_l^2}{(1 -
  x_iy_l)(1 - x_jy_l)}\rp\\ + 2\sum_{i,j}(-\theta)^{1-p(i)-q(j)}\lp
  n_\theta - m_\theta + (-\theta)^{1-p(i)} -
  (-\theta)^{1-q(j)}\rp\frac{x_iy_j}{1 - x_iy_j}.
\end{multline}
We now rewrite the first two sums in this expression such that their
dependence on the variables $x$ appear only through terms of the form
$x_iy_j/(1 - x_iy_j)$. To rewrite the first sum we use the fact that
\begin{equation*}
  \frac{x_iy_j^2}{(y_j - y_l)(1 - x_iy_j)} =
  \frac{y_l}{(y_j - y_l)(1 - x_iy_j)} + \frac{x_iy_j}{1
  - x_iy_j} - \frac{y_l}{y_j - y_l},
\end{equation*}
and the fact that the resulting sum
\begin{multline*}
  2\sum_{j\neq l}(-\theta)^{2-p(i)-q(j)-q(l)}\frac{y_l}{(y_j -
  y_l)(1 - x_iy_j)}\\ = \sum_{j\neq
  l}(-\theta)^{2-p(i)-q(j)-q(l)}\lp\frac{y_l}{(y_j - y_l)(1 - x_iy_j)}
  - \frac{y_j}{(y_j - y_l)(1 - x_iy_l)}\rp.
\end{multline*}
Setting $x = x_i$, $y = 1/y_j$ and $z = 1/y_l$ in Lemma
\ref{threeLemma} we find that
\begin{equation*}
  \frac{x_i^2y_jy_l}{(1 - x_iy_j)(1 - x_iy_l)} - \frac{y_l}{(y_j -
  y_l)(1 - x_iy_j)} + \frac{y_j}{(y_j - y_l)(1 - x_iy_l)} = 1.
\end{equation*}
We thus conclude that the first sum in \eqref{diffAction} equals
\begin{equation*}
  -2\sum_i\sum_{j\neq
   l}(-\theta)^{2-p(i)-q(j)-q(l)}\frac{x_iy_j}{1 - x_iy_j}.
\end{equation*}
The second sum in \eqref{diffAction} can be similarly rewritten to
yield
\begin{equation*}
  2\sum_{i,l}\sum_{j\neq i}(-\theta)^{2-p(i)-p(j)-q(l)}\frac{x_iy_l}{1
  - x_iy_l}.
\end{equation*}
Inserting these latter expressions for the first two sums in
\eqref{diffAction} we obtain
\begin{multline*}
  \Pi^{-1}(x,y)\lp D^2_x - \bar{D}^2_y\rp\Pi(x,y) =
  2\sum_{i,l}(-\theta)^{1-p(i)-q(l)}(n_\theta -
  m_\theta)\frac{x_iy_l}{1 - x_iy_l}\\ +
  2\sum_{i,j,l}(-\theta)^{2-p(i)-p(j)-q(l)}\frac{x_iy_l}{1 - x_iy_l} -
  2\sum_{i,j,l}(-\theta)^{2-p(i)-q(j)-q(l)}\frac{x_iy_l}{1 - x_iy_l}.
\end{multline*}
The identity \eqref{DId} for $k = 2$ now follows from the fact
that
\begin{equation*}
  \sum_j(-\theta)^{1-p(j)} = -n_\theta,\quad \sum_j(-\theta)^{1-q(j)}
  = -m_\theta.
\end{equation*}
For the remaining values of $k$ the identity \eqref{DId} can now be
obtained by first verifying that
\begin{equation*}
  D^k = \frac{1}{k+1}\left\lbrack E^0,D^{k+1}\right\rbrack,\quad
  \bar{D}^k = \frac{1}{k+1}\left\lbrack
  \bar{D}^{k+1},\bar{E}^0\right\rbrack,
\end{equation*}
and then observing that
\begin{equation*}
  \left\lbrack E^0_x,D^{k+1}_x\right\rbrack\Pi(x,y) =
  \left\lbrack\bar{D}^{k+1}_y,\bar{E}^0_y\right\rbrack\Pi(x,y).
\end{equation*}
\end{proof}

By combining the identities \eqref{EId} and \eqref{DId} we can obtain
such identities for all deformed CMS operators \eqref{CMSOps}. Indeed,
if we for each such operator $\mathcal{L}_{\bar{n}}$ let
\begin{equation*}
  \bar{\mathcal{L}}_{\bar{n},\bar{m}} =
  \sum_{k=0}^2\alpha_k\bar{D}^k_{\bar{n},\bar{m}} +
  \sum_{\ell=0}^1\beta_{\ell}\bar{E}^{\ell}_{\bar{n},\bar{m}}
\end{equation*}
then we obtain the following:

\begin{proposition}\label{IdProp}
For all deformed CMS operators $\mathcal{L}_{\bar{n}}$,
\begin{equation}\label{LId}
  \mathcal{L}_{\bar{n},x}\Pi_{\bar{n},\bar{m}}(x,y) =
  \bar{\mathcal{L}}_{\bar{n},\bar{m},y}\Pi_{\bar{n},\bar{m}}(x,y).
\end{equation}
\end{proposition}

\begin{remark}
In certain special cases these identities have been obtained
before. For $\mathcal{L}_n = D^2_n$ it is implicit in Stanley's paper
\cite{Stan}, see also Chapter VI in Macdonald \cite{Macd2}. Sergeev and
Veselov \cite{SergVes2} obtained an identity closely related to our
identity \eqref{DId} for $D^2_{\bar{n}}$ with $\bar{n}$ arbitrary and
$\bar{m} = (m,0)$ for some positive integer $m$. We also mention that
similar identities have been obtained from the point of view of
quantum Calogero-Sutherland models, see Gaudin \cite{Gaudin}, Serban
\cite{Serban} and the paper \cite{HallLang}.
\end{remark}

For the remainder of this section we set $\bar{m} = \bar{n}$. We will
show that the resulting operators
$\bar{\mathcal{L}}_{\bar{n}}:=\bar{\mathcal{L}}_{\bar{n},\bar{n}}$ are
the adjoints of $\mathcal{L}_{\bar{n}}$ for a particular scalar
product on $\Lambda_{\bar{n},\theta}$. To obtain the precise statement
we proceed in analogy with the proof of Statement 2.13 in Section VI.2
of Macdonald.

\begin{definition}
We define the scalar product $\langle\cdot,\cdot\rangle_{\bar{n}}$ on
$\Lambda_{\bar{n},\theta}$ (antilinear in its second argument) by
setting
\begin{equation}\label{sesqForm}
  \langle SP_\lambda,SP_\mu\rangle_{\bar{n}} =
  b_\lambda^{-1}\delta_{\lambda\mu}
\end{equation}
for all $\lambda,\mu\in H_{\bar{n}}$.
\end{definition}

As we now prove, for $\bar{m} = \bar{n}$ Proposition \ref{IdProp} is
equivalent to the following:

\begin{proposition}\label{dualityProposition}
For all $f,g\in\Lambda_{\bar{n},\theta}$,
\begin{equation}\label{Duality}
  \langle\mathcal{L}_{\bar{n}}f,g\rangle_{\bar{n}} = \langle
  f,\bar{\mathcal{L}}_{\bar{n}}g\rangle_{\bar{n}}.
\end{equation}
\end{proposition}

\begin{proof}
As indicated above, we prove the statement by establishing its
equivalence to Proposition \ref{IdProp} for $\bar{m} = \bar{n}$. With
\begin{equation*}
  SQ_\lambda = b_\lambda SP_\lambda
\end{equation*}
it is clear that
\begin{equation}\label{orthogonality}
  \langle SP_\lambda,SQ_\mu\rangle_{\bar{n}} =
  \delta_{\lambda\mu}
\end{equation}
for all $\lambda,\mu\in H_{\bar{n}}$. For each pair of partitions
$\lambda,\mu\in H_{\bar{n}}$ we define constants $e_{\lambda\mu}$ and
$\bar{e}_{\lambda\mu}$ by
\begin{equation*}
  \mathcal{L}_{\bar{n}} SP_\lambda = \sum_{\mu\in H_{\bar{n}}}
  e_{\lambda\mu}SP_\mu
\end{equation*}
and
\begin{equation*}
  \bar{\mathcal{L}}_{\bar{n}} SQ_\lambda = \sum_{\mu\in H_{\bar{n}}}
  \bar{e}_{\lambda\mu}SQ_\mu,
\end{equation*}
respectively. Inserting these two expressions into the identity
\eqref{LId} and using Proposition \ref{superJackExpansionProposition}
we obtain
\begin{equation*}
  \sum_{\lambda,\mu\in
  H_{\bar{n}}}e_{\lambda\mu}SP_{\mu}(x)SP_{\lambda}(y) =
  \sum_{\lambda,\mu\in
  H_{\bar{n}}}\bar{e}_{\lambda\mu}SP_{\lambda}(x)SP_{\mu}(y).
\end{equation*}
Hence, the identity \eqref{LId} is equivalent to
\begin{equation}\label{IdCondition}
  e_{\lambda\mu} = \bar{e}_{\mu\lambda}
\end{equation}
for all $\lambda,\mu\in H_{\bar{n}}$. On the other hand, using the
definition of the scalar product $\langle\cdot,\cdot\rangle_{\bar{n}}$
to compute the left- and right-hand sides of \eqref{Duality} for $f =
SP_\lambda$ and $g = SQ_\mu$ with $\lambda$ and $\mu$ running through
all partitions in $H_{\bar{n}}$, we find that also \eqref{Duality} is
equivalent to \eqref{IdCondition}. Clearly, this implies
\eqref{Duality}.
\end{proof}

It is easily verified that $\bar{E}^1_{\bar{n}} = E^1_{\bar{n}}$ and
that $\bar{D}^2_{\bar{n}} = D^2_{\bar{n}}$. As a special case of
Proposition \ref{dualityProposition} we thus obtain the following:

\begin{corollary}
Suppose that $\alpha(x) = \alpha_2x^2$ and $\beta(x) = \beta_1x$ for
some coefficients $\alpha_2$ and $\beta_1$, respectively. Then, the
corresponding deformed CMS operator $\mathcal{L}_{\bar{n}}$ is
self-adjoint for the scalar product
$\langle\cdot,\cdot\rangle_{\bar{n}}$ on $\Lambda_{\bar{n},\theta}$.
\end{corollary}

\begin{remark}
The scalar product obtained by setting $\tilde m = \tilde n = 0$ was
used by Macdonald \cite{Macd2} and the corollary can in this case be
inferred from his Statement 3.11 in Section VI.3.
\end{remark}

\section{A linear basis for the algebra $\Lambda_{\bar{n},\theta}$}
In this section we give a precise definition of the polynomials
$f^{(\bar{m})}_a$ and deduce some of their basic properties. We prove,
in particular, that a subset of these polynomials, parametrised by the
partitions in $H_{\bar{n}}\cap H_{\bar{m}}$, span the same subspace of
$\Lambda_{\bar{n},\theta}$ as the corresponding super Jack
polynomials.

The polynomials $f^{(\bar{m})}_a$ will be indexed by integer vectors
$a\in\mathbb{Z}^{|\bar{m}|}$ and before stating their definition we
introduce some notation related to this fact. We will on a number of
occasions relate such polynomials $f^{(\bar{m})}_a$ to super Jack
polynomials indexed by partitions in $H_{\bar{m}}$. For that we define
$\varphi = \varphi_{\bar{m}}$ to act on partitions $\lambda\in
H_{\bar{m}}$ as
\begin{equation*}
  \varphi(\lambda) = ({}^m\lambda,{}_m\lambda^\prime).
\end{equation*}
Note that $\varphi(\lambda)\in\mathbb{Z}^{|\bar{m}|}$ for all
$\lambda\in H_{\bar{m}}$. We let $e_i$ be the standard basis elements
in $\mathbb{Z}^{|\bar{m}|}$ defined by $(e_i)_j = \delta_{ij}$ for all
$i,j=1,\ldots,|\bar{m}|$. For each integer vector $a =
(a_1,\ldots,a_{|\bar{m}|})\in\mathbb{Z}^{|\bar{m}|}$ we let $|a| = a_1
+ \cdots + a_{|\bar{m}|}$. In addition, we will make extensive use of
the partial order on $\mathbb{Z}^{|\bar{m}|}$ defined by the
equivalence
\begin{equation}\label{reverseOrder}
  a\preceq b\Leftrightarrow a_i + \cdots + a_{|\bar{m}|}\leq b_i +
   \cdots + b_{|\bar{m}|},\quad \forall i = 1,\ldots,|\bar{m}|.
\end{equation}
By comparing definitions it is straightforward to verify that this
partial order is related to the partial order defined by
\eqref{intVectDominanceOrder} as follows:

\begin{lemma}\label{orderLemma}
Let $a,b\in\mathbb{Z}^{|\bar{m}|}$ be such that $|a| = |b|$. Then
$a\preceq b$ if and only if $b\leq a$.
\end{lemma}

We are now ready to state the following:

\begin{definition}\label{fDef}
We define the polynomials $f^{(\bar{m})}_a$,
$a\in\mathbb{Z}^{|\bar{m}|}$, by the expansion
\begin{equation}\label{faDef}
  \prod_{j<l}(1 - y_l/y_j)^{-(-\theta)^{1-q(j)-q(l)}}\prod_{i,j}(1 -
  x_iy_j)^{(-\theta)^{1-p(i)-q(j)}} = \sum_a f^{(\bar{m})}_a(x)y^a,
\end{equation}
valid for $\min_i\lp|x_i^{-1}|\rp>|y_1|>\cdots>|y_{|\bar{m}|}|$, and
where the sum is over integer vectors $a\in\mathbb{Z}^{|\bar m|}$.
\end{definition}

\begin{remark}
  These polynomials were defined by Langmann \cite{Lang} for the
  special cases $\bar{m} = \bar{n} = (N,0)$ with $N$ any positive
  integer, while the more general definition stated above is
  equivalent to the one given in \cite{HallLang}.
\end{remark}

Before studying the polynomials $f^{(\bar{m})}_a$ in general it is
instructive to consider the special cases for which $\theta = 1$ and
$\bar{m} = (M,0)$ for some positive integer $M$. In these cases they
stand in a simple relation to the so-called super Schur polynomials,
introduced by Pragacz and Thorup \cite{PragThor}; see also Fulton and
Pragacz \cite{FultPrag}.

\begin{definition}[Pragacz and Thorup]
  For each non-negative integer $k$, the polynomial $s_k(x,\tilde{x})$
  is defined by the expansion
\begin{equation*}
  \prod_{i=1}^n(1 - x_iy)^{-1}\prod_{I=1}^{\tilde{n}}(1 +
  \tilde{x}_Iy) = \sum_k s_k(x,\tilde{x})y^k.
\end{equation*}
Let $\lambda$ be a partition and let $l = \ell(\lambda)$. The
so-called super Schur polynomial $S_\lambda(x,\tilde{x})$ is then
defined as the determinant
\begin{equation}\label{superJacobiTrudi}
  S_\lambda(x,\tilde{x}) = \det\lp
  s_{\lambda_i-i+j}(x,\tilde{x})\rp_{1\leq i,j\leq l}.
\end{equation}
\end{definition}

\begin{remark}
Note that if $\tilde{n} = 0$ the polynomials $s_k$ reduce to the
complete symmetric polynomials $h_k$, the super Schur polynomials to
the `ordinary' Schur polynomials and \eqref{superJacobiTrudi} to the
so-called Jacobi-Trudi identity, originally due to Jacobi
\cite{Jacobi}; see also Section I.3 in Macdonald \cite{Macd2}.
\end{remark}

\begin{proposition}\label{superSchurProp}
  Let $\lambda$ be a partition. Suppose that $\theta = 1$ and that
  $\bar{m} = (M,0)$ for some positive integer
  $M\geq\ell(\lambda)$. Then
\begin{equation*}
  f^{(\bar{m})}_{\lambda}(x,-\tilde{x}) = S_\lambda(x,\tilde{x}).
\end{equation*}
\end{proposition}

\begin{proof}
For each $a\in\mathbb{N}^l$ we let $s_a = s_{a_1}\cdots s_{a_l}$. In
addition, we let $\delta = (l-1,l-2,\ldots,1,0)$ and observe that
\begin{equation*}
  \det\lp s_{\lambda_i-i+j}(x,\tilde{x})\rp_{1\leq i,j\leq l} =
  \det\lp s_{\lambda_i-\delta_i+\delta_j}(x,\tilde{x})\rp_{1\leq
  i,j\leq l} = \sum_{w\in S_l}\epsilon(w)s_{\lambda+\delta-w\delta},
\end{equation*}
where $S_l$ refers to the permutation group of $l$ objects and
$\epsilon(w)$ denotes the sign of the permutation $w$. Let $R_{ij}$
denote the so-called raising operator defined by $R_{ij}a = a + e_i -
e_j$ for each $a\in\mathbb{Z}^{l}$. In $\mathbb{C}\lbrack x_1^{\pm
1},\ldots,x_l^{\pm 1}\rbrack$ we then have that
\begin{equation*}
  \sum_{w\in S_l}\epsilon(w)x^{\lambda+\delta-w\delta} =
  x^{\lambda+\delta}\prod_{i<j}(x_i^{-1} - x_j^{-1}) = \prod_{i<j}(1 -
  x_i/x_j)x^\lambda = \prod_{i<j}(1 - R_{ij})x^\lambda.
\end{equation*}
By applying the linear map $\mathbb{C}\lbrack x_1^{\pm
  1},\ldots,x_l^{\pm 1}\rbrack\rightarrow \Lambda_{\bar{n},\theta}$
defined by $x^a\mapsto s_a$ for all $a\in\mathbb{Z}^l$ we find that
\begin{equation*}
  S_\lambda = \prod_{i<j}(1 - R_{ij})s_\lambda.
\end{equation*}
On the other hand,
\begin{equation*}
\begin{split}
  \sum_{a\in\mathbb{Z}^M}f^{(\bar{m})}(x,-\tilde{x})y^a &=
  \prod_{j<l}(1 - y_l/y_j)\prod_{i,j}(1 - x_iy_j)^{-1}\prod_{I,j}(1 +
  \tilde{x}_Iy_j)\\ & = \prod_{j<l}(1 -
  y_l/y_j)\sum_{a\in\mathbb{N}^M}s_a(x,\tilde{x})y^a.
\end{split}
\end{equation*}
Comparing coefficients of $y^\lambda$ we obtain
\begin{equation*}
  f^{(\bar{m})}_{\lambda}(x,-\tilde{x}) = \prod_{i<j}(1 -
  R_{ij})s_\lambda(x,\tilde{x}) = S_\lambda(x,\tilde{x}).
\end{equation*}
\end{proof}

We proceed to prove that the set of polynomials
$f^{(\bar{m})}_{\varphi(\lambda)}$ and the set of super Jack
polynomials $SP_\lambda$, both indexed by the partitions $\lambda\in
H_{\bar{n}}\cap H_{\bar{m}}$ of a given weight $|\lambda|$, span the
same linear subspace of $\Lambda_{\bar{n},\theta}$. The main facts
required for the proof are contained in the following:

\begin{lemma}\label{transitionLemma}
  For any non-negative integer $k$ there exists a transition matrix $M
  = (M_{a\mu})$, defined by the equalities
\begin{equation*}
  f^{(\bar{m})}_a = \sum_\mu M_{a\mu}SP_\mu,
\end{equation*}
from the polynomials $\left(f^{(\bar{m})}_a\right)$, indexed by the
integer vectors $a\in\mathbb{Z}^{|\bar m|}$ such that $|a| = k$, to
the super Jack polynomials $(SP_\mu)$, indexed by the partitions
$\mu\in H_{\bar{n}}\cap H_{\bar{m}}$ of weight $|\mu| = k$. The
entries $M_{a\mu}$ of $M$ are non-zero only if $\varphi(\mu)\preceq
a$. In addition, if $a = \varphi(\lambda)$ for some partition
$\lambda\in H_{\bar{n}}\cap H_{\bar{m}}$,
\begin{equation*}
  M_{a\lambda} = (-1)^{|{}_m\lambda|}
  b_\lambda(\theta)b_{({}_m\lambda^\prime)}\left(\theta^{-1}\right).
\end{equation*}
\end{lemma}

\begin{proof}
For any Laurent polynomial $g = g(y)$ in the variables $y =
(y_1,\ldots,y_{|\bar{m}|})$ and integer vector $a\in\mathbb{Z}^{|\bar
m|}$ we let $\lbrack g\rbrack_a$ denote the coefficient of $y^a$ in
$g(y)$. It follows from the definition of the polynomials
$f^{(\bar{m})}_a$ and Proposition \ref{superJackExpansionProposition}
that
\begin{equation*}
  f^{(\bar{m})}_a(x) = \sum_{\mu\in H_{\bar{n}}\cap
  H_{\bar{m}}}b_\mu(\theta)\left\lbrack\prod_{j<l}(1 -
  y_l/y_j)^{-(-\theta)^{1-q(j)-q(l)}}SP_\mu(y)\right\rbrack_a
  SP_\mu(x).
\end{equation*}
Expanding the products in binomial series we conclude that the entries
\begin{equation*}
  M_{a\mu} = b_\mu(\theta)\prod_{j<l}\sum_{\nu_{jl}=0}^\infty
  (-1)^{\nu_{jl}}\binom{-(-\theta)^{1-q(j)-q(l)}}{\nu_{jl}}\left\lbrack
  SP_\mu(y) \right\rbrack_{a + \sum_{i<j}\nu_{jl}(e_j-e_l)}.
\end{equation*}
We recall from Lemma \ref{homLemma} that a super Jack polynomial
$SP_\mu$, if non-zero, is homogeneous of degree $|\mu|$. Lemmas
\ref{triangLemma} and \ref{orderLemma} thus imply that such a super
Jack polynomial $SP_\mu(y)$ is a linear combination of monomials $x^b$
with $b\succeq\varphi(\mu)$. Since
\begin{equation*}
  a + \sum_{j<l}\nu_{jl}(e_j - e_l)\preceq a
\end{equation*}
for all non-negative integers $\nu_{ij}$, it follows that the entries
$M_{a\mu}$ are non-zero only if $|\mu| = |a|$ and $\varphi(\mu)\preceq
a$. Clearly, this establishes both the existence, as well as the
stated triangular structure, of the transition matrix $M$. Suppose
that the integer vector $a = \varphi(\lambda)$ for some partition
$\lambda\in H_{\bar{n}}\cap H_{\bar{m}}$. It is then clear that the
only term which contributes to $M_{a\lambda}$ is the leading term
$\lbrack SP_\lambda(y)\rbrack_{\varphi(\lambda)}$ of $SP_\lambda$. The
statement thus follows from Lemma \ref{triangLemma}.
\end{proof}

\begin{corollary}\label{faCorollary}
A polynomial $f^{(\bar{m})}_a$, $a\in\mathbb{Z}^{|\bar{m}|}$, is
non-zero only if $a\succeq 0$. In that case,
$f^{(\bar{m})}_a\in\Lambda_{\bar{n},\theta}$ and it is homogeneous of
degree $|a|$.
\end{corollary}

\begin{proof}
  We observe that if $a$ violates the condition $a\succeq 0$ there
  exist no partition $\mu\in H_{\bar{m}}$ such that
  $\varphi(\mu)\preceq a$. The statement is thus a direct consequence
  of Lemmas \ref{homLemma} and \ref{transitionLemma}.
\end{proof}

Using Lemma \ref{transitionLemma} we now prove the following:

\begin{proposition}\label{spanningProposition}
Assume that $\theta$ is not a negative rational number or zero. Then,
as $\lambda$ runs through all partitions of a given weight $|\lambda|
= k$ in $H_{\bar{n}}\cap H_{\bar{m}}$, the corresponding polynomials
$f^{(\bar{m})}_{\varphi(\lambda)}$ form a linear basis for the linear
space
\begin{equation*}
  \mathbb{C}\left\langle SP_\lambda : \lambda\in H_{\bar{n}}\cap
  H_{\bar{m}}, |\lambda| =
  k\right\rangle\subseteq\Lambda^k_{\bar{n},\theta}.
\end{equation*}
In particular, if $m\geq n$ and $\tilde m\geq \tilde n$, then they
form a linear basis for $\Lambda^k_{\bar{n},\theta}$.
\end{proposition}

\begin{proof}
  We let $K = (K_{\lambda\mu})$ be the restriction of the transition
  matrix $M$ defined in Lemma \ref{transitionLemma} to the polynomials
  $\left(f^{(\bar{m})}_a\right)$, indexed by the integer vectors $a =
  \varphi(\lambda)$ for some partition $\lambda\in H_{\bar{n}}\cap
  H_{\bar{m}}$ of weight $|\lambda| = k$. We note that $K$ is a square
  matrix, whereas $M$ is not. We also note that, in the terminology of
  Section I.6 in Macdonald \cite{Macd2}, $K$ is a strictly upper
  triangular matrix in the sense that $K_{\lambda\mu} = 0$ unless
  $\varphi(\lambda)\succeq\varphi(\mu$). It follows from Lemma
  \ref{transitionLemma} and the definition of $b_\lambda(\theta)$ that
  all diagonal elements of $K$ (under the stated condition on
  $\theta$) are all well defined and non-zero. It is readily verified
  that the inverse of such a matrix exists and is of the same form,
  c.f.\ Statement 6.1 in Macdonald (\textit{loc.\ cit.}). For all
  $\lambda\in H_{\bar{n}}\cap H_{\bar{m}}$ we thus have
\begin{equation*}
  SP_\lambda = \sum_\mu
  (K^{-1})_{\lambda\mu}f^{(\bar{m})}_{\varphi(\mu)}
\end{equation*}
and the statement follows.
\end{proof}

\section{Polynomial eigenfunctions of deformed CMS operators}
For each of the deformed CMS operators \eqref{CMSOps} we construct in
this section sets of polynomial eigenfunctions, parametrised by the
partitions in $H_{\bar{n}}\cap H_{\bar{m}}$, and expressed as linear
combinations of the polynomials $f_a^{(\bar{m})}$. Under a certain
condition of non-degeneracy on their eigenvalues we prove that these
polynomial eigenfunctions span the same subspace of
$\Lambda_{\bar{n},\theta}$ as the super Jack polynomials labeled by
the partitions in $H_{\bar{n}}\cap H_{\bar{m}}$.

The first step in our construction is to compute the action of the
operators $\mathcal{L}_{\bar{n}}$ on the polynomials
$f^{(\bar{m})}_a$. Before proceeding to the computation we sketch our
approach. We observe that the generating function for the polynomials
$f^{(\bar{m})}_a$, as defined by \eqref{faDef}, is related to the
function $\Pi_{\bar{n},\bar{m}}$ as follows:
\begin{equation*}
  \prod_{j<l}(1 - y_l/y_j)^{-(-\theta)^{1-q(j)-q(l)}}\prod_{i,j}(1 -
  x_iy_j)^{(-\theta)^{1-p(i)-q(j)}} =
  G_{\bar{m}}(y)\Pi_{\bar{n},\bar{m}}(x,y)
\end{equation*}
with
\begin{equation*}
  G_{\bar{m}}(y) = \prod_{j<l}(1 -
  y_l/y_j)^{-(-\theta)^{1-q(j)-q(l)}}.
\end{equation*}
We also observe that Proposition \ref{IdProp} implies the following
identity:
\begin{equation*}
  \mathcal{L}_{\bar{n},x}G_{\bar{m}}(y)\Pi_{\bar{n},\bar{m}}(x,y) =
  \left(G_{\bar{m}}\bar{\mathcal{L}}_{\bar{n},\bar{m}}
  G^{-1}_{\bar{m}}\right)_y G_{\bar{m}}(y)\Pi_{\bar{n},\bar{m}}(x,y).
\end{equation*}
By first expanding the function
$G_{\bar{m}}(y)\Pi_{\bar{n},\bar{m}}(x,y)$ as in \eqref{faDef}, then
computing the right-hand side of this identity, and finally comparing
coefficients with the left-hand side, we can thus obtain the action of
the operators $\mathcal{L}_{\bar{n}}$ on the polynomials
$f^{(\bar{m})}_a$. The difficult part of this computation will be to
obtain the explicit form of the operators
$G_{\bar{m}}\bar{\mathcal{L}}_{\bar{n},\bar{m}}G_{\bar{m}}^{-1}$. However,
the computation is considerably simplified if we first consider only
the operators $E^{\ell}_{\bar{n}}$ and $D^k_{\bar{n}}$.

Before proceeding we mention that a detailed discussion of
the approach sketched above in the case of the `ordinary' CMS
operators $\mathcal{L}_n$, as well as a comparison with a construction of
their eigenfunctions in terms of more elementary bases for the
symmetric polynomials, can be found in \cite{Hallnas}.

We find it convenient to introduce the following notation: to each
$a\in\mathbb{Z}^{|\bar{m}|}$ we associate the `shifted' integer vector
$a^+ = a + s$ with the shift $s = (s_1,\ldots,s_{|\bar{m}|})$ given by
\begin{equation*}
  s_j = (-\theta)^{1-q(j)}\left((-\theta)^{-q(j)}(j - 1) - n\right) -
  (-\theta)^{-q(j)}(m + \tilde{n}) + m.
\end{equation*}
In addition, for each non-negative integer $\nu$ and $k = 0,1,2$ we
let
\begin{equation*}
  E^{k\nu}_{jl} = (2 - k + \nu)e_l -\nu e_j.
\end{equation*}
We are now ready to prove the following:

\begin{lemma}\label{actionLemma}
Let $a\in\mathbb{Z}^{|\bar{m}|}$ be such that $a\succeq 0$. Then
\begin{subequations}\label{EAction}
\begin{align}
  E^1_{\bar{n}}f^{(\bar{m})}_a &= |a|f^{(\bar{m})}_a,\\
  E^0_{\bar{n}}f^{(\bar{m})}_a &= \sum_{j=1}^{|\bar{m}|}(a^+_j -
  1)f^{(\bar{m})}_{a-e_j}.
\end{align}
\end{subequations}
Furthermore, with $k = 0,1$,
\begin{subequations}\label{DAction}
\begin{align}
  D^2_{\bar{n}}f^{(\bar{m})}_a &=
  \left(\sum_{j=1}^{|\bar{m}|}(-\theta)^{q(j)}a_j(a_j - 1 + 2s_j) +
  2|a|\right)f^{(\bar{m})}_a\\ \nonumber &\quad + 2(\theta -
  1)\sum_{j<l}(-\theta)^{1-q(j)-q(l)}\sum_{\nu=1}^\infty \nu
  f^{(\bar{m})}_{a-E^{2\nu}_{jl}},\\ D^k_{\bar{n}}f^{(\bar{m})}_a & =
  \sum_j(-\theta)^{q(j)}\left(a^+_j - 2\delta_{k0} +
  k\left((-\theta)^{-q(j)} - 1\right)\right)(a^+_j -
  1)f^{(\bar{m})}_{a-(2-k)e_j}\\ \nonumber &\quad + (\theta -
  1)\sum_{j<l}(-\theta)^{1-q(j)-q(l)}\sum_{\nu=0}^\infty(2\nu + 2 -
  k)f^{(\bar{m})}_{a-E^{k\nu}_{jl}}.
\end{align}
\end{subequations}
\end{lemma}

\begin{proof}
We will throughout the proof assume that, whenever they occur, $\ell =
0,1$ and $k = 0,1,2$. For simplicity of notation we let $\upsilon_j =
(-\theta)^{-q(j)}$ for all $j = 1,\ldots,|\bar{m}|$ and suppress the
index's $\bar{n}$ and $\bar{m}$. In addition, we let $M_G$ denote the
operator of multiplication by the function $G$.

In order to prove \eqref{EAction} we note that the similarity
transformation of a partial derivative $\partial/\partial y_j$ by
$M_G$ is given by
\begin{equation*}
  M_G\frac{\partial}{\partial y_j}M^{-1}_G = \frac{\partial}{\partial
  y_j} - \theta\sum_{l>j}\upsilon_j\upsilon_l\frac{y_l/y_j}{y_j - y_l}
  - \theta\sum_{l<j}\upsilon_j\upsilon_l\frac{1}{y_j - y_l}.
\end{equation*}
Multiplying by $y_j^{2-\ell}$, taking the sum over $j$, and
interchanging summation index's in the resulting second sum, we obtain
\begin{equation*}
\begin{split}
  M_G \bar{E}^\ell_yM^{-1}_G &=
  \sum_jy_j^{2-\ell}\frac{\partial}{\partial y_j} + (1 -
  \ell)\sum_j\lp n_\theta -
  \theta\sum_{l<j}\upsilon_l\rp\upsilon_jy_j\\ &=
  \sum_jy_j^{2-\ell}\frac{\partial}{\partial y_j} + (1 -
  \ell)\sum_js_jy_j.
\end{split}
\end{equation*}
It follows from Lemma \ref{IdLemma} and the definition of the
polynomials $f_a$ that
\begin{equation*}
  \sum_a\lp E^\ell_xf_a(x)\rp y^a = \sum_a
  f_a(x)G(y)\bar{E}^\ell_yG^{-1}(y)y^a.
\end{equation*}
By first computing the right-hand side and then comparing coefficients
of the monomials $y^a$ with the left-hand side we thus obtain
\eqref{EAction}.

We proceed to prove \eqref{DAction}. Using the identity
\begin{equation}\label{aSimpleId}
  \frac{y_j^{n-1}y_l}{y_j-y_l} = \frac{y_j^n}{y_j - y_l} - y_j^{n-1},
\end{equation}
valid for any integer $n$, it is readily verified that the similarity
transformation of an operator $\bar{D}^k$ by $M_G$ is given by
\begin{equation}\label{simTrans}
\begin{split}
  M_G\bar{D}^k_yM^{-1}_G &=
  \sum_j\upsilon_j^{-1}y_j^{4-k}\frac{\partial^2}{\partial y_j^2} +
  \sum_j\left(C_{j,k} +
  2\theta\sum_{l>j}\upsilon_l\right)y_j^{3-k}\frac{\partial}{\partial
  y_j}\\ &\quad + G(y)\bar{D}^k_yG^{-1}(y).
\end{split}
\end{equation}
Repeatedly applying the identity \eqref{aSimpleId}, and using the fact
that $s_j\upsilon_j^{-1} = n_\theta - \theta\sum_{l<j}\upsilon_l$, we
deduce through a straightforward but somewhat lengthy computation
\begin{multline*}
  G(y)\bar{D}^k_y G^{-1}(y) = \theta\sum_{j\neq l}\lp k\upsilon_l(1 -
  \upsilon_j)\frac{y_j^{2-k}y_l}{y_j - y_l} -
  (\theta\upsilon_j\upsilon_l - 1)\upsilon_l\frac{y_j^{2-k}y_l^2}{(y_j
  - y_l)^2}\rp\\ - \theta^2\sum_{\substack{j\neq l,l^\prime\\ l\neq
  l^\prime}}\upsilon_j\upsilon_l\upsilon_{l^\prime}
  \frac{y_j^{2-k}y_ly_{l^\prime}}{(y_j - y_l)(y_j - y_{l^\prime})} -
  2\theta n_\theta\sum_{j\neq l}\upsilon_j\upsilon_l
  \frac{y_j^{2-k}y_l}{y_j - y_l}\\ + \sum_j(s_j(s_j + 1) -
  n_\theta\upsilon_j(n_\theta\upsilon_j +1) + k(\upsilon_j - 1)(s_j -
  n_\theta\upsilon_j))\upsilon^{-1}_j y_j^{2-k} + P_k(y).
\end{multline*}
Since $2(\upsilon_lf - \upsilon_jg) = (\upsilon_l +
\upsilon_j)(f + g) + (\upsilon_l - \upsilon_j)(f - g)$ for any two
functions $f$ and $g$, we have
\begin{equation*}
\begin{split}
  \theta\sum_{j\neq l}(\theta\upsilon_j\upsilon_l -
  1)\upsilon_l\frac{y_j^{2-k}y_l^2}{(y_j - y_l)^2} &=
  \frac{\theta}{2}\sum_{j<l}(\theta\upsilon_j\upsilon_l -
  1)(\upsilon_l + \upsilon_j)\frac{y_j^{2-k}y_l^2 +
  y_j^2y_l^{2-k}}{(y_j - y_l)^2}\\ &\quad +
  \frac{\theta}{2}\sum_{j<l}(\theta\upsilon_j\upsilon_l -
  1)(\upsilon_l - \upsilon_j)\frac{y_j^{2-k}y_l^2 -
  y_j^2y_l^{2-k}}{(y_j -y_l)^2}.
\end{split}
\end{equation*}
By a straightforward computation using the fact that $2\upsilon_l(1 -
\upsilon_j) = \upsilon_l(1 - \upsilon_j)(1 - \theta\upsilon_j)$ we
thus obtain
\begin{multline*}
  \theta\sum_{j\neq l}\lp k\upsilon_l(1 -
  \upsilon_j)\frac{y_j^{2-k}y_l}{y_j - y_l} -
  (\theta\upsilon_j\upsilon_l - 1)\upsilon_l\frac{y_j^{2-k}y_l^2}{(y_j
  - y_l)^2}\rp\\ =
  -\frac{\theta}{2}\sum_{j<l}(\theta\upsilon_j\upsilon_l -
  1)(\upsilon_l + \upsilon_j)\frac{y_j^{2-k}y_l^2 +
  y_j^2y_l^{2-k}}{(y_j - y_l)^2}\\ - \delta_{k,2}\theta\sum_{j<
  l}\lp\frac{1}{2}(\theta\upsilon_j\upsilon_l + 1)(\upsilon_j +
  \upsilon_l) - (1 + \theta)\upsilon_j\upsilon_l\rp.
\end{multline*}
Rewriting the sum such that its summand is symmetric in the summation
index's and setting $x = 1/y_j$, $y = 1/y_l$ and $z = 1/y_{l^\prime}$
in Lemma \ref{threeLemma} we deduce
\begin{equation*}
  \theta^2\sum_{\substack{j\neq l,l^\prime\\l\neq
  l^\prime}}\upsilon_j\upsilon_l\upsilon_{l^\prime}
  \frac{y_j^{2-k}y_ly_{l^\prime}}{(y_j - y_l)(y_j - y_{l^\prime})} =
  \delta_{k,2}2\theta^2
  \sum_{j<l<l^\prime}\upsilon_j\upsilon_l\upsilon_{l^\prime}.
\end{equation*}
Using the fact that
\begin{equation*}
  2\sum_{j\neq l}\upsilon_j\upsilon_l\frac{y_j^{2-k}y_l}{y_j - y_l} =
\sum_{j\neq l}\upsilon_j\upsilon_l(\delta_{k,0}y_jy_l - \delta_{k,2}),
\end{equation*}
as well as the definition of the polynomials $P_k(y)$, we thus obtain
\begin{equation}\label{conjugation}
\begin{split}
  G(y)\bar{D}^k_yG^{-1}(y) &= (1 - \delta_{k2})\sum_j\big(s_j(s_j + 1)+
  k(\upsilon_j - 1)s_j\big)\upsilon_j^{-1}y_j^{2-k}\\ &\quad -
  \frac{\theta}{2}\sum_{j<l}(\theta\upsilon_j\upsilon_l -
  1)(\upsilon_j + \upsilon_l)\left(\frac{y_j^{2-k}y_l^2 +
  y_j^2y_l^{2-k}}{(y_j - y_l)^2} - \delta_{k2}\right).
\end{split}
\end{equation}
We now make the Laurent series expansion
\begin{equation*}
  \frac{y_j^{2-k}y_l^2 + y_j^2y_l^{2-k}}{(y_j - y_l)^2} - \delta_{k2}
  = \sum_{\nu=0}^\infty (2\nu + 2 - k)\frac{y_l^{\nu+2-k}}{y_j^\nu},
\end{equation*}
valid for
$\min_i\lp|x_i^{-1}|\rp>|y_1|>\cdots>|y_{\bar{m}}|$. Inserting this
Laurent series expansion, as well as \eqref{conjugation}, into
\eqref{simTrans} and using the fact that
\begin{equation*}
  C_{j,k} + 2\theta\sum_{l>j}\upsilon_l = \big(2(s_j + 1) +
  k(\upsilon_j - 1)\big)\upsilon_j^{-1}
\end{equation*}
it is straightforward to compute $G(y)\bar{D}^k_yG^{-1}(y)y^a$ for all
integer vectors $a\in\mathbb{Z}^{|\bar{m}|}$ such that $a\succeq
0$. It is then readily verified that \eqref{DAction} follows from
Lemma \ref{IdLemma}, as in the proof of \eqref{EAction} above.
\end{proof}

By Corollary \ref{faCorollary} and Proposition
\ref{spanningProposition} we thus obtain the following:

\begin{corollary}\label{degreeCorollary}
Set $\Lambda^N_{\bar{n},\theta} = \lbrace 0\rbrace$ for all
$N<0$. Then, for $\ell = 0,1$,
\begin{equation*}
  E^\ell_{\bar{n}} : \Lambda^N_{\bar{n},\theta}\rightarrow
  \Lambda^{N-(1-\ell)}_{\bar{n},\theta}
\end{equation*}
and, for $k = 0,1,2$,
\begin{equation*}
  D^k_{\bar{n}} : \Lambda^N_{\bar{n},\theta}\rightarrow
  \Lambda^{N-(2-k)}_{\bar{n},\theta}
\end{equation*}
for all non-negative integers $N$.
\end{corollary}

\begin{remark}
We say that a linear operator $L$ on $\Lambda_{\bar{n},\theta}$ is
homogeneous of degree $k$ if
$L\Lambda^N_{\bar{n},\theta}\subset\Lambda^{N-k}_{\bar{n},\theta}$
for all $N\geq 0$. Corollary \ref{degreeCorollary} then states that
each operator $E^\ell_{\bar{n}}$ and $D^k_{\bar{n}}$ is homogeneous
of degree $\ell - 1$ and $k - 2$, respectively. In addition, using
Proposition \ref{dualityProposition} it is readily inferred that
$\bar{E}^{\ell}_{\bar{n},\bar{m}}$ and $\bar{D}^k_{\bar{n},\bar{m}}$
are similarly homogeneous of degree $3-\ell$ and $4-k$,
respectively.
\end{remark}

As mentioned above, the polynomial eigenfunctions we construct for the
deformed CMS operators \eqref{CMSOps} will be labeled by the
partitions in $H_{\bar{n}}\cap H_{\bar{m}}$. For a given such
partition $\lambda$ the corresponding eigenfunction will be obtained
as a linear combination of polynomials $f^{(\bar{m})}_b$ indexed by
integer vectors of the form $b = \varphi(\lambda)-a\succeq 0$, where
the integer vectors $a$ will be contained in certain subsets of
$\lbrace a\in\mathbb{Z}^{|\bar{m}|}: a\succeq 0\rbrace$. Although it
is difficult to give a precise characterisation of these subsets, it
will become evident below that they are contained in sets of integer
vectors $\mathscr{C}_{|\bar{m}|}(\alpha,\beta)$ which are determined
in a simple manner by the two polynomials $\alpha$ and $\beta$, and
which in many cases are significantly smaller than $\lbrace
a\in\mathbb{Z}^{|\bar{m}|}: a\succeq 0\rbrace$.

\begin{definition}\label{coneDef}
Fix a non-negative integer $n$. Associate to each pair of polynomials
$\alpha$ and $\beta$ the index sets $I_\alpha =
\lbrace k: \alpha_k\neq 0\rbrace$ and $I_\beta = \lbrace \ell:
\beta_\ell\neq 0\rbrace$. For $k = 0,1,2$ let
\begin{equation*}
  A^{(k)}_n = \lbrace a\in\mathbb{Z}^n: a\succeq 0, |a| = (2 -
  k)\nu~\text{for some}~\nu\in\mathbb{N}\rbrace,
\end{equation*}
and for $\ell = 0,1$ let
\begin{equation*}
  B^{(\ell)}_n = \lbrace a\in\mathbb{N}^n: |a| = (1 -
  \ell)\nu~\text{for some}~\nu\in\mathbb{N}\rbrace.
\end{equation*}
We then define $\mathscr{C}_n = \mathscr{C}_n(\alpha,\beta)$ to be the
set of integer vectors of the form
\begin{equation*}
  a = \sum_{k\in I_\alpha}a^{(k)} + \sum_{\ell\in I_\beta}b^{(\ell)}
\end{equation*}
for some $a^{(k)}\in A^{(k)}_n$ and $b^{(\ell)}\in B^{(\ell)}_n$.
\end{definition}

\begin{remark}
We note that in most special cases this definition can be
considerably simplified, e.g., for $\alpha(x) = x^2$ and $\beta =
0$, corresponding to the deformed CMS operator
$\mathcal{L}_{\bar{n}} = \bar{D}^2_{\bar{n}}$, we have
$\mathscr{C}_n = A^{(2)}_n$. Hence, in this special case,
$\mathscr{C}_n$ is the set of all integer vectors $a\succeq 0$ such
that $|a| = 0$.
\end{remark}

At this point we fix the two polynomials $\alpha$ and $\beta$ and
consider the resulting deformed CMS operator $\mathcal{L}_{\bar{n}}$,
as defined by \eqref{CMSOps}. It is clear from Lemma \ref{actionLemma}
and Corollary \ref{faCorollary} that $\mathcal{L}_{\bar{n}}$ maps each
polynomial $f^{(\bar{m})}_a$ to a linear combination of polynomials
$f^{(\bar{m})}_{a-b}$ with $b\in \mathscr{C}_{|\bar{m}|}$ such that
$b\preceq a$. This suggests that for each partition $\lambda\in
H_{\bar{n}}\cap H_{\bar{m}}$ there exists an eigenfunction
$P^{(\bar{m})}_\lambda$ of $\mathcal{L}_{\bar{n}}$ of the form
\begin{equation}\label{eigenFunc}
  P^{(\bar{m})}_\lambda = \sum_a u_\lambda(a)
  f^{(\bar{m})}_{\varphi(\lambda)-a},
\end{equation}
where the sum is over integer vectors $a\in\mathscr{C}_{|\bar{m}|}$
such that $0\preceq a\preceq\varphi(\lambda)$. Indeed, using Lemma
\ref{actionLemma}, it is straightforward to verify that this is the
case if the coefficients $u_\lambda(a)$ satisfy the recurrence
relation
\begin{multline}\label{recursRels}
  \big(\mathcal{E}_{\bar{m}}(\varphi(\lambda)) -
  \mathcal{E}_{\bar{m}}(\varphi(\lambda)-a)\big) u_\lambda(a) =
  \alpha_0\sum_j(-\theta)^{q(j)}a^+_j\left(a^+_j + 1\right)u_\lambda(a
  + 2e_j)\\ + \sum_ja^+_j\left(\beta_0 +
  \alpha_1(-\theta)^{q(j)}\left(a^+_j +
  (-\theta)^{-q(j)}\right)\right)u_\lambda(a + e_j)\\ + (\theta -
  1)\sum_{j<l}(-\theta)^{1-q(j)-q(l)}
  \sum_{k=0}^2\alpha_k\sum_{\nu=0}^\infty (2\nu + 2 -
  k)u_\lambda\left(a + E^{k\nu}_{jl}\right)
\end{multline}
with
\begin{equation*}
  \mathcal{E}_{\bar{m}}(a) =
  \alpha_2\sum_{j=1}^{|\bar{m}|}(-\theta)^{q(j)}a_j(a_j - 1 + 2s_j) +
  (2\alpha_2 + \beta_1)|a|.
\end{equation*}
In addition, the corresponding eigenvalues are given by
$\mathcal{E}_{\bar{m}}(\varphi(\lambda))$. We note that if
$\mathcal{E}_{\bar{m}}(\varphi(\lambda))\neq
\mathcal{E}_{\bar{m}}(\varphi(\lambda)-a)$ for all
$a\in\mathscr{C}_{|\bar{m}|}$ such that $0\prec
a\preceq\varphi(\lambda)$ then the recursions relation
\eqref{recursRels} uniquely determines the coefficients $u_\lambda(a)$
once the leading coefficient $u_\lambda(0)$ has been fixed. We also
note that since $0\preceq a\preceq\varphi(\lambda)$, the sum over
$\nu$ in \eqref{recursRels} truncates after a finite number of terms.

\begin{definition}\label{admissibleDef}
We say that a partition $\lambda\in H_{\bar{m}}$ is
$\bar{m}$-admissible if $\mathcal{E}(\varphi(\lambda))\neq
\mathcal{E}(\varphi(\lambda) - a)$ for all
$a\in\mathscr{C}_{|\bar{m}|}$ such that $0\prec
a\preceq\varphi(\lambda)$.
\end{definition}

A polynomial $P^{(\bar{m})}_\lambda$ with the leading coefficient
$u_\lambda(0)$ fixed to some non-zero constant and remaining
coefficients determined by the recursion relation \eqref{recursRels}
is thus a well-defined eigenfunction of the deformed CMS operator
$\mathcal{L}_{\bar{n}}$ if the partition $\lambda$ is
$\bar{m}$-admissible. An important aspect of our construction is that
these eigenfunctions can be normalised such that they are independent
of $\bar{m}$. To give a first indication of this fact we proceed to
show that, although the expression for $\mathcal{E}_{\bar{m}}(a)$
above depends on $\bar{m}$, the eigenvalues
$\mathcal{E}_{\bar{m}}(\varphi(\lambda))$ of the eigenfunctions
$P^{(\bar{m})}_\lambda$ do not. For the convenience of the reader we
provide a complete proof of this fact although a proof has been
previously given in \cite{HallLang}; see Lemma 4.1.

\begin{lemma}
For all partitions $\lambda\in H_{\bar{n}}\cap H_{\bar{m}}$,
\begin{equation*}
\begin{split}
  \mathcal{E}_{\bar{m}}(\varphi(\lambda)) &=
  \mathcal{E}_{(\ell(\lambda),0)}(\lambda)\\ &=
  \alpha_2\sum_j\lambda_j\big(\lambda_j + 1 + 2\theta(n - j + 1) -
  2\tilde{n}\big) + \beta_1\sum_j\lambda_j,
\end{split}
\end{equation*}
where the two sums run over all parts of $\lambda$.
\end{lemma}

\begin{proof}
By the definition of the map $\varphi$ we have
\begin{equation*}
  \mathcal{E}(\varphi(\lambda)) - (2\alpha_2+\beta_1)|\lambda| =
  \alpha_2\sum_{j=1}^m\lambda_j(\lambda_j - 1 + 2s_j) -
  \alpha_2\theta\sum_{J=1}^{\tilde{m}}
  {}_m\lambda_J^\prime({}_m\lambda_J^\prime - 1 + 2s_{m+J}).
\end{equation*}
For any partition $\mu$ it is easily verified that
\begin{equation*}
  \sum_i i\mu_i^\prime = \frac{1}{2}\sum_i\mu_i(\mu + 1),
\end{equation*}
where the two sums run over the parts of $\mu^\prime$ and $\mu$,
respectively. The definition of the `shift' $s$ thus imply
\begin{equation*}
  \theta\sum_{J=1}^{\tilde{m}}{}_m\lambda_J^\prime
  \big({}_m\lambda_J^\prime - 1 + 2s_{m+J}) =
  -\sum_J{}_m\lambda_J({}_m\lambda_J - 1 + 2\theta(n - m - J + 1) -
  2\tilde{n}\big)
\end{equation*}
where the latter sum is over the parts of ${}_m\lambda$. We now obtain
the statement by using again the definition of the shift $s$.
\end{proof}

For the remainder of this section we set the leading coefficients
\begin{equation}\label{leadingCoeff}
  u_\lambda(0) = (-1)^{|{}_m\lambda|}
  b^{-1}_{{}_m\lambda^\prime}\left(\theta^{-1}\right)
\end{equation}
for all polynomials $P^{(\bar{m})}_\lambda$ with $\lambda\in
H_{\bar{n}}\cap H_{\bar{m}}$ an $\bar{m}$-admissible partition. This
is motivated by the following:

\begin{proposition}\label{differentValuesProp}
Let $\bar{k}\in\mathbb{N}^2$ and let $\lambda\in
H_{\bar{n}}\cap H_{\bar{m}}\cap H_{\bar{k}}$ be $\bar{m}$- and
$\bar{k}$-admissible. Then $P^{(\bar{k})}_\lambda =
P^{(\bar{m})}_\lambda$.
\end{proposition}

\begin{proof}
To prove the statement we will make use of the fact that any super
Jack polynomial
\begin{equation}\label{superJackExp}
  SP_\mu = \sum_a c_{\mu a}f^{(\bar{m})}_a
\end{equation}
for some coefficients $c_{\mu a}$ and where the sum is over integer
vectors $a\in\mathbb{Z}^{|\bar{m}|}$ such that $|a| = |\mu|$ and
$a\preceq\varphi(\mu)$. This fact can be inferred either from Lemma
\ref{transitionLemma} (c.f.\ the proof of Proposition
\ref{spanningProposition}) or directly from Theorem \ref{superJackThm}
below. We let $C_{\bar{m}}(\mu)$ be the set of integer vectors
$a\prec\varphi_{\bar{m}}(\mu)$ such that $\varphi_{\bar{m}}(\mu) -
a\in\mathscr{C}_{|\bar{m}|}$. From \eqref{superJackExp} and Lemma
\ref{actionLemma} we deduce that
\begin{equation}\label{prelActionOnSuperJacks}
  \mathcal{L}_{\bar{n}}SP_\mu = \sum_{|a|=|\mu|}\mathcal{E}(a)c_{\mu
  a}f^{(\bar{m})}_a + \sum_{b\in C_{\bar{m}}(\mu)}c^\prime_{\mu
  b}f^{(\bar{m})}_b
\end{equation}
for some coefficients $c^\prime_{\mu b}$. To proceed we consider the
two cases $\alpha_2 = 0$ and $\alpha_2\neq 0$ separately. We first
assume that $\alpha_2 = 0$. In this case $\mathcal{E}_{\bar{m}}(a) =
\mathcal{E}_{\bar{m}}(\varphi(\mu))$ for all integer vectors
$a\in\mathbb{Z}^{|\bar{m}|}$ such that $|a| = |\mu|$. Hence,
\begin{equation*}
  \mathcal{L}_{\bar{n}}SP_\mu = \mathcal{E}_{\bar{m}}(\varphi(\mu))SP_\mu
  + \sum_{b\in C_{\bar{m}}(\mu)}c^\prime_{\mu
    b}f^{(\bar{m})}_b.
\end{equation*}
It is easily inferred from Definition \ref{coneDef} that if $b\in
C_{\bar{m}}$ then any integer vector $a\in\mathbb{Z}^{|\bar{m}|}$ such
that $|a| = |b|$ and $a\preceq b$ is contained in $C_{\bar{m}}$. It
thus follows from Lemma \ref{transitionLemma} that
\begin{equation}\label{actionOnSuperJacks}
  \mathcal{L}_{\bar{n}}SP_\mu = \mathcal{E}(\varphi(\mu))SP_\mu +
  \sum_{\nu\in C_{\bar{m}}(\mu)} d_{\mu\nu}SP_\nu
\end{equation}
for some coefficients $d_{\mu\nu}$. We turn now to the case
$\alpha_2\neq 0$. In this case all integer vectors $a\neq
\varphi(\mu)$ in \eqref{prelActionOnSuperJacks} are contained in
$C_{\bar{m}}$. It follows that
\begin{equation*}
  \mathcal{L}_{\bar{n}}SP_\mu = \mathcal{E}_{\bar{m}}(\varphi(\mu))SP_\mu 
  + \sum_{b\in C_{\bar{m}}(\mu)}c^{\prime\prime}_{\mu  b}f^{(\bar{m})}_b
\end{equation*}
for some coefficients $c^{\prime\prime}_{\mu b}$. Referring again to
Lemma \ref{transitionLemma} we thus find that
\eqref{actionOnSuperJacks} holds true also for $\alpha_2\neq 0$. It is
clear that if we replace $\bar{m}$ by $\bar{k}$ in the discussion
above we obtain \eqref{actionOnSuperJacks} with $C_{\bar{m}}(\mu)$
replaced by $C_{\bar{k}}(\mu)$. Since $\lambda$ is both $\bar{m}$- and
$\bar{k}$-admissible, it follows that with $C(\lambda)$ either of the
sets $C_{\bar{m}}(\lambda)$, $C_{\bar{k}}(\lambda)$ or
$C_{\bar{m}}(\lambda)\cap C_{\bar{k}}(\lambda)$ there exist a unique
eigenfunction $P_\lambda$ of $\mathcal{L}_{\bar{n}}$ such that
\begin{equation}\label{jackEigenFuncs}
  P_\lambda = SP_\lambda + \sum_{\mu\in C(\lambda)}
  u_{\lambda\mu}SP_\mu
\end{equation}
for some coefficients $u_{\lambda\mu}$. Since
$C_{\bar{m}}(\lambda)\cap C_{\bar{k}}(\lambda)$ is contained in both
$C_{\bar{m}}(\lambda)$ and $C_{\bar{k}}(\lambda)$, these eigenfunctions
must all coincide. It thus follows from \eqref{eigenFunc},
\eqref{leadingCoeff} and Lemma \ref{transitionLemma} that
\begin{equation*}
  P^{(\bar{k})}_\lambda = P_\lambda^{(\bar{m})} =
  b_\lambda(\theta)P_\lambda.
\end{equation*}
\end{proof}

Given a partition $\lambda\in H_{\bar{n}}$ we thus obtain different
series representations for the same eigenfunction of
$\mathcal{L}_{\bar{n}}$ by varying the value of $\bar{m}$ such that
$\lambda\in H_{\bar{m}}$ is $\bar{m}$-admissible. It is interesting to
note that the complexity of the resulting series representation is in
many cases highly dependent on the specific value we choose for
$\bar{m}$. To give a simple and concrete example of this fact we
consider the deformed CMS operator $\mathcal{L}_{\bar{n}} =
D^2_{\bar{n}}$ for $\bar{n} = (2,1)$. Suppose that we are interested
the eigenfunction indexed by the partition $\lambda = (1^3)\equiv
(1,1,1)$. We observe that $\mathscr{C}_3(x^2,0) = A_3^{(2)}$ consists
of all integer vectors $a\in\mathbb{Z}^3$ such that $a\succeq 0$ and
$|a| = 0$. In the discussion below we assume that the parameter
$\theta$ is such that the partition $\lambda$ is $\bar{m}$-admissible
for all values of $\bar{m}$ under consideration. That such values of
$\theta$ exist follows from Proposition \ref{admissibleProp} below. If
we set $\bar{m} = \bar{n} = (2,1)$ then we obtain a series
representation of the form
\begin{equation*}
\begin{split}
  P^{(2,1)}_{(1^3)} &= u_{(1^3)}(0)f^{(2,1)}_{(1^3)} +
  u_{(1^3)}((-1,0,1))f^{(2,1)}_{(2,1)} +
  u_{(1^3)}((-1,1))f^{(2,1)}_{(2,0,1)}\\ &\quad +
  u_{(1^3)}((-2,2))f^{(2,1)}_{(3,-1,1)} +
  u_{(1^3)}((-2,1,1))f^{(2,1)}_{(3)}.
\end{split}
\end{equation*}
On the other hand, with $\bar{m} = (0,1)$ the same eigenfunction
is given by the series
\begin{equation*}
  P^{(0,1)}_{(1^3)} = u_{(1^3)}(0)f^{(0,1)}_{(3)},
\end{equation*}
containing only the polynomial
\begin{equation*}
  f^{(0,1)}_{(3)} = -\binom{-1/\theta}{1}x_1x_2\tilde{x}_1 -
  \binom{-1/\theta}{2}(x_1 + x_2)\tilde{x}_1^2 -
  \binom{-1/\theta}{3}\tilde{x}_1^2.
\end{equation*}
In this case there is clearly a large difference in the complexity of
the series representation obtained for $\bar{m} = (2,1)$ and $(0,1)$,
respectively. In general, the 'simplest' series representation for a
given eigenfunction is obtained by choosing $\bar{m}$ such that
$|\bar{m}|$ is minimised. This reflects the fact that the complexity
of an eigenfunction depends to a large extent on the partition to
which it corresponds, and to a lesser extent on the number of
variables $\bar{n}$. For a further discussion of the complexity of
this type of series representations for the eigenfunctions of
(deformed) CMS operators, and their dependence on the value of
$\bar{m}$, see \cite{HallLang}, in particular the discussion following
Theorem 4.1.

In analogy with Proposition \ref{spanningProposition} we proceed to
prove that the eigenfunctions $P^{(\bar{m})}_\lambda$, indexed by the
partitions $\lambda\in H_{\bar{n}}\cap H_{\bar{m}}$, span the same
linear subspace of $\Lambda_{\bar{n},\theta}$ as the corresponding
super Jack polynomials.

\begin{theorem}\label{basisThm}
  Assume that $\theta$ is not a negative rational number or zero and
  that all partitions in $H_{\bar{n}}\cap H_{\bar{m}}$ are
  $\bar{m}$-admissible. Then, the corresponding eigenfunctions
  $P^{(\bar{m})}_\lambda$ of the deformed CMS operator
  $\mathcal{L}_{\bar{n}}$ form a linear basis for the linear space
\begin{equation*}
  \mathbb{C}\langle SP_\lambda: \lambda\in H_{\bar{n}}\cap
  H_{\bar{m}}\rangle\subseteq\Lambda_{\bar{n},\theta}.
\end{equation*}
In particular, if $m\geq n$ and $\tilde{m}\geq\tilde{n}$, then they
form a linear basis for $\Lambda_{\bar{n},\theta}$.
\end{theorem}

\begin{proof}
Fix a partition $\lambda\in H_{\bar{n}}\cap H_{\bar{m}}$ and
consider the eigenfunctions $\left(P^{(\bar{m})}_\mu\right)$,
indexed by the partitions $\mu\in H_{\bar{n}}\cap H_{\bar{m}}$ such
that $\varphi(\mu)\preceq\varphi(\lambda)$. Observe that
$\varphi(\lambda) - a\preceq \varphi(\lambda)$ for all
$a\in\mathscr{C}_{|\bar{m}|}$. It follows from \eqref{eigenFunc} and
Lemma \ref{transitionLemma} that there exist a well defined strictly
upper triangular transition matrix $N = (N_{\lambda\mu})$ (c.f. the
proof of Proposition \ref{spanningProposition}) from the
eigenfunctions $\left(P^{(\bar{m})}_\mu\right)$ to the super Jack
polynomials $(SP_\mu)$, also indexed by the partitions $\mu\in
H_{\bar{n}}\cap H_{\bar{m}}$ such that
$\varphi(\mu)\preceq\varphi(\lambda)$. Furthermore, the specific
form chosen for the leading coefficients $u_\lambda(0)$ together
with Lemma \ref{transitionLemma} and the definition of
$b_\lambda(\theta)$ implies that all of its diagonal elements are
non-zero (under the stated assumption on $\theta$). Hence, its
inverse $N^{-1}$ exist and is of the same form. For each partition
$\mu\in H_{\bar{n}}\cap H_{\bar{m}}$ such that
$\varphi(\mu)\preceq\varphi(\lambda)$ we thus have that
\begin{equation*}
  SP_\mu = \sum_\nu(N^{-1})_{\mu\nu}P^{(\bar{m})}_\nu.
\end{equation*}
Since $\lambda$ can be fixed to any partition in $H_{\bar{n}}\cap
H_{\bar{m}}$ the statement follows.
\end{proof}

We conclude this section by establishing two sufficient conditions for
all partitions in $H_{\bar{n}}\cap H_{\bar{m}}$ to be
$\bar{m}$-admissible.

\begin{proposition}\label{admissibleProp}
Assume that at least one of the coefficients $\alpha_2$ and $\beta_1$
are non-zero. All partitions in $H_{\bar{n}}\cap H_{\bar{m}}$ are then
$\bar{m}$-admissible if either of the following two conditions
are satisfied:
\begin{enumerate}
\item the coefficient $\alpha_2$ is zero,\label{cond1}
\item $\theta$ is not a negative rational number,
$\tilde{m}\leq 1$, and the coefficients $\alpha_1$, $\alpha_0$ and
$\beta_0$ are zero.\label{cond2}
\end{enumerate}
\end{proposition}

\begin{proof}
For any partition $\lambda\in H_{\bar{n}}\cap H_{\bar{m}}$ and integer
vector $a\in\mathscr{C}_{\bar{m}}$ such that $0\prec
a\preceq\varphi(\lambda)$,
\begin{multline*}
  \mathcal{E}_{\bar{m}}(\varphi(\lambda)) -
  \mathcal{E}_{\bar{m}}(\varphi(\lambda) - a) =
  2\alpha_2\sum_{j=1}^{|\bar{m}|}(-\theta)^{q(j)}a_j
  \big((\varphi(\lambda))_j + s_j\big)\\ -
  \alpha_2\sum_{j=1}^{|\bar{m}|}(-\theta)^{q(j)}a_j(a_j + 1) +
  (2\alpha_2 + \beta_1)|a|.
\end{multline*}
Suppose that $\alpha_2 = 0$. Then $|a|>0$ and Condition \eqref{cond1}
is clearly sufficient for all partitions in $H_{\bar{n}}\cap
H_{\bar{m}}$ to be $\bar{m}$-admissible. Suppose instead that
$\alpha_2\neq 0$. Assume furthermore that $\tilde{m}\leq 1$ and that
$\alpha_1$, $\alpha_0$ and $\beta_0$ are zero. Then it is readily
verified that the integer vector
\begin{equation*}
  a = -\sum_{j<l}\nu_{jl}(e_j - e_l)
\end{equation*}
for some non-negative integers $\nu_{jl}$ not all zero. We thus have
\begin{equation*}
\begin{split}
  \sum_{j=1}^{|\bar{m}|}(-\theta)^{q(j)}a_j\big((\varphi(\lambda))_j +
  s_j\big) &= -\sum_{1\leq j<l\leq m}\nu_{jl}(\lambda_j + s_j -
  \lambda_l - s_l)\\ &\quad- \sum_{j=1}^m\nu_{jm+1}\big(\lambda_j +
  s_j + \theta((\varphi(\lambda))_{m+1} + s_{m+1})\big).
\end{split}
\end{equation*}
In addition,
\begin{equation*}
  \sum_{j=1}^{|\bar{m}|}(-\theta)^{q(j)}a_j(a_j + 1) =
  \sum_{j=1}^ma_j(a_j + 1) - \theta a_{m+1}(a_{m+1} + 1).
\end{equation*}
The definition of the `shift' $s$ and the fact that $a_{m+1} =
\sum_{j=1}^m\nu_{jm+1}$ thus imply
\begin{equation*}
  \mathcal{E}_{\bar{m}}(\varphi(\lambda))
  - \mathcal{E}_{\bar{m}}(\varphi(\lambda) - a) =
  -K(\lambda,a) - \theta L(\lambda,a)
\end{equation*}
with
\begin{equation*}
  K(a,\lambda) = \sum_{j=1}^ma_j(a_j + 1) + 2\sum_{1\leq j<l\leq
  m}\nu_{jl}(\lambda_j - \lambda_l) + 2\sum_{j=1}^m\nu_{jm+1}\lambda_j
\end{equation*}
and
\begin{equation*}
  L(a,\lambda) = 2\sum_{1\leq j<l\leq m}\nu_{jl}(l - j) +
  2\sum_{j=1}^m\nu_{jm+1}\left((\varphi(\lambda))_{m+1} -
  \frac{1}{2}a_{m+1} + m - j\right).
\end{equation*}
Since $a\preceq\varphi(\lambda)$, $(\varphi(\lambda))_{m+1} -
a_{m+1}\geq 0$. It follows that both $K(a,\lambda)$ and $L(a,\lambda)$
are strictly positive integers. Clearly, this implies that also
Condition \eqref{cond2} is sufficient for all partitions in
$H_{\bar{n}}\cap H_{\bar{m}}$ to be $\bar{m}$-admissible.
\end{proof}

\section{Series expansions of the super Jack polynomials}
A recursion relation closely related to \eqref{recursRels} is in
\cite{HallLang} solved by iteration. Rather than deriving a similar
solution of \eqref{recursRels}, we concentrate in this section on the
special case $\mathcal{L}_{\bar{n}} = D^2_{\bar{n}}$, corresponding to
$\alpha(x) = x^2$ and $\beta(x) = 0$. In particular, we exhibit the
explicit solution to the recursion relation \eqref{recursRels} for
this special case and prove that the resulting eigenfunctions
\eqref{eigenFunc}, up to a constant of proportionality, coincide with
the super-Jack polynomials. We note that the super Jack polynomials
are known to be eigenfunctions of the differential operator
$D^2_{\bar{n}}$; see e.g.\ Sergeev and Veselov
\cite{SergVes2}. However, since the eigenvalues of $D^2_{\bar{n}}$ in
general are not distinct, this fact alone does not imply that the
eigenfunctions \eqref{eigenFunc} we have constructed are proportional
to the super-Jack polynomials.

\begin{theorem}\label{superJackThm}
Let $\lambda\in H_{\bar{n}}\cap H_{\bar{m}}$ be an
$\bar{m}$-admissible partition. Then
\begin{equation}\label{superJackExpansion}
\begin{split}
  b_\lambda SP_\lambda &= f^{(\bar{m})}_{\varphi(\lambda)} +
  \sum_{s=1}^\infty 2^s(\theta -
  1)^s\sum_{j_1<l_1}(-\theta)^{1-q(j_1)-q(l_1)}
  \sum_{\nu_1=1}^\infty\nu_1\\ &\quad
  \times\cdots\times\sum_{j_s<l_s}(-\theta)^{1-q(j_s)-q(l_s)}
  \sum_{\nu_s=1}^\infty\nu_s\\
  &\quad\times\prod_{r=1}^s\lp\mathcal{E}_{\bar{m}}(\varphi(\lambda)) -
  \mathcal{E}_{\bar{m}}\lp\varphi(\lambda) -
  \sum_{t=r}^s\nu_tE_{j_tl_t}\rp\rp^{-1}
  f^{(\bar{m})}_{\varphi(\lambda)-\sum_{r=1}^s\nu_rE_{j_rl_r}}
\end{split}
\end{equation}
with
\begin{equation*}
  E_{jl} = e_l - e_j.
\end{equation*}
\end{theorem}

\begin{proof}
Setting $\alpha_2 = 1$, and the remaining coefficients $\alpha_k$ and
$\beta_\ell$ to zero, in \eqref{recursRels} we find that
\eqref{eigenFunc} is an eigenfunction of the differential operator
$D^2_{\bar{n}}$ if the coefficients $u_\lambda(a)$ satisfy the
recursion relation
\begin{equation}\label{recursRelD2}
  \left(\mathcal{E}_{\bar{m}}(\varphi(\lambda)) -
  \mathcal{E}_{\bar{m}}(a)\right)u_\lambda(a) = 2(\theta -
  1)\sum_{j<l}(-\theta)^{1-q(j)-q(l)}\sum_{\nu=1}^\infty\nu
  u_\lambda(a+\nu(e_l - e_j)).
\end{equation}
We set $u_\lambda(0) = 1$ and write $\delta_a(b)$ for the Kronecker
delta of two integer vectors $a,b\in\mathbb{Z}^{|\bar{m}|}$, i.e.,
$\delta_a(b)$ equals $1$ if $a = b$ and zero otherwise. By iterating
\eqref{recursRelD2} we thus deduce that each coefficient
\begin{equation*}
\begin{split}
  u_\lambda(a) &= \delta_{\varphi(\lambda)}(a) +
  \sum_{s=1}^\infty 2^s(\theta - 1)^s\sum_{j_1<l_1}\sum_{\nu_1=1}^\infty
  \frac{\nu_1(-\theta)^{1-q(j_1)-q(l_1)}}{\mathcal{E}_{\bar{m}}(\varphi(\lambda))
  - \mathcal{E}_{\bar{m}}(a)}\\ &\quad \times\sum_{j_2<l_2}\sum_{\nu_2=1}^\infty
  \frac{\nu_2(-\theta)^{1-q(j_2)-q(l_2)}}{\mathcal{E}_{\bar{m}}(\varphi(\lambda))
  - \mathcal{E}_{\bar{m}}(a + \nu_1E_{j_1l_1})}\times\cdots\\
  &\quad\times\sum_{j_s<l_s} \sum_{\nu_s=1}^\infty
  \frac{\nu_s(-\theta)^{1-q(j_s)-q(l_s)}}{\mathcal{E}_{\bar{m}}(\varphi(\lambda))
  - \mathcal{E}_{\bar{m}}(a+\sum_{r=1}^{s-1}\nu_rE_{j_rl_r})}
  \delta_{\varphi(\lambda)}\lp a + \sum_{r=1}^s\nu_rE_{j_rl_r}\rp\\ &=
  \delta_{\varphi(\lambda)}(a) + \sum_{s=1}^\infty 2^s(\theta -
  1)^s\sum_{j_1<l_1}(-\theta)^{1-q(j_1)-q(l_1)}
  \sum_{\nu_1=1}^\infty\nu_1\\
  &\quad\times\cdots\times\sum_{j_s<l_s}(-\theta)^{1-q(j_s)-q(l_s)}
  \sum_{\nu_s=1}^\infty\nu_s\\
  &\quad\times\prod_{r=1}^s\lp\mathcal{E}_{\bar{m}}(\varphi(\lambda)) -
  \mathcal{E}_{\bar{m}}\lp\varphi(\lambda) -
  \sum_{t=r}^s\nu_tE_{j_tl_t}\rp\rp^{-1}\delta_{\varphi(\lambda)}\lp a
  + \sum_{r=1}^s\nu_rE_{j_rl_r}\rp.
\end{split}
\end{equation*}
Inserting this into \eqref{eigenFunc} we obtain the right-hand side of
\eqref{superJackExpansion}. It remains to show that the eigenfunction
equals $b_\lambda SP_\lambda$. It follows from Lemma
\ref{transitionLemma} that the right-hand side of
\eqref{superJackExpansion} is a linear combination of super Jack
polynomials $SP_\mu$ with $\varphi(\lambda) -
\varphi(\mu)\in\mathscr{C}_{|\bar{m}|}$ and leading term $b_\lambda
SP_\lambda$. Moreover, the assumption that $\lambda$ is
$\bar{m}$-admissible implies that all terms other than the leading one
have eigenvalues different from
$\mathcal{E}_{\bar{m}}(\varphi(\lambda))$. Consequently, these terms
must be zero.
\end{proof}

\begin{remark}
It is important to note that Corollary \ref{faCorollary} implies that
the series expansion \eqref{superJackExpansion} of the super Jack
polynomials only contains a finite number of non-zero terms and thus
is well defined.
\end{remark}

We proceed to further study a number of particularly simple special
cases of the series expansion \eqref{superJackExp} of the super Jack
polynomials. In particular, it is interesting to consider this series
expansion for $\theta = 1$. In that case $b_\lambda = 1$, and the
right hand side of \eqref{superJackExp} contains only the leading term
$f^{(\bar{m})}_{\varphi(\lambda)}$. Setting $\bar{m} = (m,0)$, and
using Proposition \ref{superSchurProp}, as well as case (2) of
Proposition \ref{admissibleProp}, we thus recover the following:

\begin{corollary}
For $\theta = 1$ and each partition $\lambda\in H_{\bar{n}}$,
\begin{equation*}
	SP_\lambda(x,-\tilde{x}) = S_\lambda(x,\tilde{x}).
\end{equation*}
\end{corollary}

In addition, allowing arbitrary integer vectors $\bar{m}$ we obtain
the following generalisation of Proposition \ref{superSchurProp}:

\begin{corollary}\label{faSchurCor}
For $\theta = 1$ and each $\bar{m}$-admissible partition $\lambda\in
H_{\bar{n}}\cap H_{\bar{m}}$,
\begin{equation*}
  f^{(\bar{m})}_{\varphi(\lambda)}(x,-\tilde{x}) = S_\lambda(x,\tilde{x}).
\end{equation*}
\end{corollary}

We note that, in general, the series expansion \eqref{superJackExp}
contains not only polynomials $f^{(\bar{m})}_a$ parametrised by an
integer vector of the form $a = \varphi(\lambda)$ for some partition
$\lambda\in H_{\bar{n}}\cap H_{\bar{m}}$, i.e., the sum is not only
over the polynomials singled out in Proposition
\ref{spanningProposition} for providing a linear basis for a natural
subspace of $\Lambda_{\bar{n},\theta}$. However, if $|\bar{m}| \leq 2$
it follows from Corollary \ref{faCorollary} that the sum in fact run
only over precisely the $f^{(\bar{m})}_a$ corresponding to a partition
in $H_{\bar{n}}\cap H_{\bar{m}}$ under the map $\varphi$. In
addition, for $|\bar{m}| \leq 2$ the polynomials $f^{(\bar{m})}_a$
themselves are particularly simple. More precisely, for $|\bar{m} =
1|$ we have the following:

\begin{proposition}\label{m1Prop}
Suppose that $|\bar{m}| = 1$. Let $\lambda\in H_{\bar{n}}\cap
H_{\bar{m}}$. Then
\begin{equation*}
  b_\lambda SP_\lambda = f^{(\bar{m})}_{\varphi(\lambda)}.
\end{equation*}
Moreover, for each (positive) integer $a$,
\begin{align*}
  f^{(1,0)}_{(a)}(x,\tilde{x};\theta) &= \sum_{r=0}^a (-1)^r e_r(\tilde{x})g_{a-r}(x;\theta),\\
  f^{(0,1)}_{(a)}(x,\tilde{x};\theta) &= \sum_{r=0}^a (-1)^r
  e_r(x)g_{a-r}(\tilde{x};1/\theta).
\end{align*}
\end{proposition}

\begin{proof}
We first observe that $\mathcal{C}_{\bar{m}}$ contains only the zero
vector. Consequently, all partitions $\lambda\in H_{\bar{m}}$ are
trivially $\bar{m}$-admissible. The first part of the statement is
thus immediate from Theorem \ref{superJackThm}. To prove the second
part of the statement we first note that, by definition,
\begin{equation*}
  \frac{\prod_I(1-\tilde{x}_Iy_1)}{\prod_i(1-x_iy_1)^\theta}
  = \sum_{a\geq 0}f^{(1,0)}_{(a)}(x,\tilde{x};\theta)y_1^a.
\end{equation*}
On the other hand,
\begin{equation*}
  \frac{\prod_I(1-\tilde{x}_Iy_1)}{\prod_i(1-x_iy_1)^\theta}
  = \sum_{r\geq 0}(-1)^r e_r(\tilde{x})y_1^r\sum_{s\geq 0}g_s(x;\theta)y_1^s.
\end{equation*}
The formula for the polynomials $f^{(1,0)}_{(a)}$ is now obtained by
comparing coefficients in the two expansions above. The formula for
the polynomials $f^{(0,1)}_{(a)}$ can be verified in a similar manner.
\end{proof}

Furthermore, for $|\bar{m}| = 2$ the following statement holds true:

\begin{proposition}\label{m2Prop}
Suppose that $|\bar{m}| = 2$. Let $\lambda\in H_{\bar{n}}\cap H_{\bar{m}}$ be an
$\bar{m}$-admissible partition. Then
\begin{equation*}
\begin{split}
  b_{\lambda}SP_{\lambda} &= f^{(\bar{m})}_{\varphi(\lambda)} +
  \sum_\nu 2^{\ell(\nu)}(\theta -
  1)^{\ell(\nu)}(-\theta)^{\frac{1}{2}(m-\tilde{m})\ell(\nu)}\nu_1\cdots
  \nu_{\ell(\nu)}\\ &\quad
  \times\prod_{r=1}^{\ell(\nu)}\lp\mathcal{E}_{\bar{m}}(\varphi(\lambda))
  - \mathcal{E}_{\bar{m}}\lp\varphi(\lambda) -
  |\nu|_rE_{12}\rp\rp^{-1}
  f^{(\bar{m})}_{\varphi(\lambda)-|\nu|E_{12}},
\end{split}
\end{equation*}
where the sum is over all positive integer vectors $\nu =
(\nu_1,\ldots,\nu_{\ell(\nu)})$ such that $|\nu| \leq
(\varphi(\lambda))_2$, and where we have used the notation
\begin{equation*}
	|\nu|_r =\nu_r + \cdots + \nu_s
\end{equation*}
for all positive integers $r\leq\ell(\nu)$. Moreover, for each
(positive) integer vector $a = (a_1,a_2)$,
\begin{equation*}
  f^{(\bar{m})}_a = \sum_{t=0}^{a_2}\binom{-(-\theta)^{\frac{1}{2}(m-\tilde{m})}}{t}
  p^{(\bar{m})}_{(a_1+t,a_2-t)}
\end{equation*}
with
\begin{align*}
  p^{(2,0)}_a &= \sum_{r=0}^{a_1}(-1)^r e_r(\tilde{x})g_{a_1-r}(x;\theta)\sum_{s=0}^{a_2}(-1)^s e_s(\tilde{x})g_{a_2-s}(x;\theta),\\
  p^{(1,1)}_a &= \sum_{r=0}^{a_1}(-1)^r e_r(\tilde{x})g_{a_1-r}(x;\theta)\sum_{s=0}^{a_2}(-1)^s e_s(x)g_{a_2-s}(\tilde{x};1/\theta),\\
  p^{(0,2)}_a &= \sum_{r=0}^{a_1}(-1)^r
  e_r(x)g_{a_1-r}(\tilde{x};1/\theta)\sum_{s=0}^{a_2}(-1)^s
  e_s(x)g_{a_2-s}(\tilde{x};1/\theta).
\end{align*}
\end{proposition}

\begin{proof}
The first part of the statement is just a reformulation of Theorem
\ref{superJackThm}. To prove the second part of the statement we
first consider the case $\bar{m} = (2,0)$. The corresponding
polynomials $f^{(\bar{m})}_a$ are then defined by the expansion
\begin{equation*}
  (1 - y_2/y_1)^{-\theta}\frac{\prod_I(1-\tilde{x}_Iy_1)(1-\tilde{x}_Iy_2)}{\prod_i(1-x_iy_1)^{\theta}(1-x_iy_2)^{\theta}} = \sum_a f^{\bar{m}}_a(x,\tilde{x})y^a.
\end{equation*}
On the other hand,
\begin{multline*}
  \frac{\prod_I(1-\tilde{x}_Iy_1)(1-\tilde{x}_Iy_2)}{\prod_i(1-x_iy_1)^{\theta}(1-x_iy_2)^{\theta}}\\
  = \sum_{r\geq 0}(-1)^r e_r(\tilde{x})y_1^r \sum_{s\geq 0}(-1)^s
  e_s(\tilde{x})y_2^s \sum_{p\geq 0}g_p(x;\theta)y_1^p \sum_{q\geq
    0}g_q(x;\theta)y_2^q.
\end{multline*}
The formula for the polynomials $f^{(2,0)}_a$ is now
obtained by expanding the factor $(1 - y_2/y_1)^{-\theta}$ and
comparing coefficients in the two expansions above. The two remaining
cases $\bar{m} = (1,1)$ and $(0,2)$ are proved similarly.
\end{proof}

Finally, we note that the expansions just obtained are non-trivial and
interesting already in the special case $\theta = 1$. For example,
from Corollary \ref{faSchurCor} and Proposition \ref{m2Prop} we deduce
the following expansion for the super Schur polynomials labeled by
hook partitions:
\begin{equation*}
  S_{(a_1,1^{a_2})}(x,\tilde{x}) = \sum_{t=0}^{a_2}(-1)^{a_2-t}\left(\sum_{r=0}^{a_1+t}e_r(\tilde{x})h_{a_1+t-r}(x)\sum_{s=0}^{a_2-t}e_s(x)h_{a_2-t-s}(\tilde{x})\right),
\end{equation*}
where we have used the fact that both the elementary symmetric
polynomials $e_r$, as well as the complete symmetric polynomials
$h_r$, are homogeneous of degree $r$. For the 'ordinary' Schur
polynomials, i.e., for $\tilde{x} = 0$ this reduces to the well known
formula
\begin{equation*}
  s_{(a_1,1^{a_2})}(x) = \sum_{t=0}^{a_2}(-1)^{a_2-t}e_{a_2-t}(x)h_{a_1+t}(x);
\end{equation*}
see e.g.\ Example 9 in Section I.3 of Macdonald \cite{Macd2}.

\section{Concluding remarks}
We have in this paper studied the polynomial eigenfunctions of the
deformed CMS operators \eqref{CMSOps} in terms of series expansion in
the polynomials $f^{(\bar{m})}_a$. In particular, we have obtained
(under a certain condition of non-degeneracy on their eigenvalues) a
linear basis for these eigenfunctions. In addition, we have
demonstrated in the special case of the super Jack polynomials that
these series expansions in the polynomials $f^{(\bar{m})}_a$ are
rather explicit. In conclusion, we briefly discuss two related papers
and some remaining problems.

As mentioned in the introduction, the present paper is closely related
to a recent paper by Langmann and the author \cite{HallLang}. In this
latter paper, eigenfunctions of both `ordinary' CMS- as well as
deformed CMS operators are obtained from the point of view of quantum
many-body systems of Calogero-Sutherland type. Many of the results
obtained here are results of questions raised in this paper. We also
mention a recent paper by Langmann \cite{Lang3}, which offers an
alternative interpretation of his original construction of the
polynomial eigenfunctions of the operator $\mathcal{L}_n$
corresponding to the so-called Sutherland model.

In Section 6 we showed that the eigenfunctions $P^{(\bar{m})}_\lambda$
can be normalised such that they are independent of the specific value
of $\bar{m}$. For the polynomials $f^{(\bar{m})}_a$ the situation is
more complicated. In fact, it is largely an open problem to
characterise their dependence on $\bar{m}$. However, certain
properties can be inferred from the results of the present paper. For
example, it follows from Lemma \ref{transitionLemma} that the
polynomials $f^{(\bar{m})}_{\varphi(\lambda)}$, corresponding to a
given partition $\lambda$, have the same leading term when expanded in
super Jack polynomials. Also, it is easily inferred from Definition
\ref{fDef} that with a fixed integer vector $a =
(a_1,\ldots,a_{|\bar{m}|})$ the corresponding polynomial
$f^{(\bar{m})}_a$ is invariant under the replacement of $\bar{m} =
(m,\tilde{m})$ by $(m,\tilde{m}+N)$ for any positive integer $N$. We
note, however, that it is easily verified in special cases that the
polynomials $f^{\bar{m}}_a$ can be distinctly different for
different values of $\bar{m}$; c.f.\ the discussion following
Proposition \ref{differentValuesProp}.

As noted in the discussion preceding Proposition \ref{m1Prop}, we
have in many instances worked with an 'overcomplete' set of
polynomials $f^{(\bar{m})}_a$. It would be desirable to be able to
rewrite the corresponding expression such that the only involve a set
of linearly independent polynomials $f^{(\bar{m})}_a$, e.g.\ those
singled out in Proposition \ref{spanningProposition}. In fact, already
Lemma \ref{transitionLemma} provides enough information to obtain the
structure of the expansion of an arbitrary polynomial
$f^{(\bar{m})}_a$ in these latter polynomials. However, it provides
very little insight into the explicit form of the coefficients in such
an expansion.

Finally, we observe that in many statements in Sections 6 and 7 it is
required that a certain partition be $\bar{m}$-admissible. This
condition encodes the level of non-degeneracy required of the
eigenvalues $\mathcal{E}(\lambda)$. In Proposition
\ref{admissibleProp} we proved two sufficient conditions for all
partitions in $H_{\bar{n}}\cap H_{\bar{m}}$ to be
$\bar{m}$-admissible. In particular, for $\alpha_2 = 0$ this is always
the case. On the other hand, for $\alpha_2\neq 0$ a number of special
cases remain to be investigated. An important such case, corresponding
to the super Jack polynomials, is that for which only $\alpha_2$ is
non-zero and $\bar{m}$ is arbitrary.
\newline

\emph{Acknowledgments.} I would like to thank E. Langmann for a
number of helpful discussions and for his comments on a preliminary
version of the paper. Financial support from the European Union
through the FP6 Marie Curie RTN ENIGMA (Contract number
MRTN-CT-200405652) is also gratefully acknowledged.

\bibliographystyle{amsalpha}

\providecommand{\bysame}{\leavevmode\hbox to3em{\hrulefill}\thinspace}
\providecommand{\href}[2]{#2}

\end{document}